\documentclass[letterpaper]{article} 
\usepackage{authblk}
\usepackage{arxiv}  
\usepackage{times}  
\usepackage{helvet}  
\usepackage{courier}  
\usepackage[hyphens]{url}  
\usepackage{graphicx} 
\usepackage{natbib}
\usepackage{ulem}
\usepackage[table]{xcolor}
\usepackage{color-edits}
\addauthor[Hui]{hui}{blue}
\usepackage{float} 
\usepackage{mathtools}
\usepackage{bm}
\usepackage{wrapfig}
\usepackage{booktabs}
\usepackage{tikz}
\usepackage[utf8]{inputenc} 
\usepackage[T1]{fontenc}    
\usepackage{booktabs}       
\usepackage{amsfonts}       
\usepackage{nicefrac}       
\usepackage{microtype}      
\usepackage{lipsum}         

\usepackage{subcaption}
\usepackage{amsmath}
\usepackage{algorithm}
\usepackage{algorithmic}
\usepackage{pgfgantt}
\usepackage{appendix}
\usepackage{mathrsfs}
\usepackage{multirow}

\usepackage{adjustbox}
\usepackage{threeparttable}
\usepackage{hyperref}       
\usepackage{tabularx,cleveref}
\usepackage{amsthm}
\usepackage{pifont}
\usepackage{enumitem}
\makeatletter
\newsavebox{\@tabnotebox}

\makeatother

\title{D-PDLP: Scaling PDLP to Distributed Multi-GPU Systems}

\theoremstyle{definition}

\usepackage{siunitx}
\theoremstyle{remark}
\author{%
  \textbf{Hongpei Li$^{1}$\thanks{Work completed during an internship at Cardinal Operations.}}, Yicheng Huang$^2$, Huikang Liu$^3$, Dongdong Ge$^3$, Yinyu Ye$^4$\\
  $^1$Cardinal Operations \\ 
  $^2$Shanghai University of Finance and Economics \\ 
  $^3$Shanghai Jiao Tong University \\ 
  $^4$Stanford University \\
  \texttt{ishongpeili@gmail.com, hkl1u@sjtu.edu.cn, ddge@sjtu.edu.cn}
}

\AtBeginDocument{%
  \addtolength\abovedisplayskip{-0.05\baselineskip}%
  \addtolength\belowdisplayskip{-0.1\baselineskip}%
  \addtolength\abovedisplayshortskip{-0.05\baselineskip}%
  \addtolength\belowdisplayshortskip{-0.1\baselineskip}%
}

\usepackage{titlesec}\titlespacing*{\section}{0pt}{-0.05\baselineskip}{-0.05\baselineskip}
\titlespacing*{\subsection}{0pt}{-0.075\baselineskip}{-0.075\baselineskip}
\titlespacing*{\subsubsection}{0pt}{-0.025\baselineskip}{-0.05\baselineskip}
\begin{document}

\maketitle

\begin{abstract}
We present a distributed framework of the Primal-Dual Hybrid Gradient (PDHG) algorithm for solving massive-scale linear programming (LP) problems. Although PDHG-based solvers demonstrate strong performance on single-node GPU architectures, their applicability to industrial-scale instances is often limited by single-GPU computational throughput. To overcome these challenges, we propose \textbf{D-PDLP}, the first \underline{D}istributed \underline{PDLP} framework, which extends PDHG to a multi-GPU setting via a practical two-dimensional grid partitioning of the constraint matrix. To improve load balance and computational efficiency, we introduce a block-wise random permutation strategy combined with nonzero-aware matrix partitioning. By distributing the intensive computation required in PDHG iterations, the proposed framework harnesses multi-GPU parallelism to achieve substantial speedups with relatively low communication overhead. Extensive experiments on standard LP benchmarks (including MIPLIB and Mittelmann instances) as well as huge-scale real-world datasets show that our distributed implementation, built upon cuPDLPx~\citep{lu2024cupdlpx}, achieves strong scalability and high performance while preserving full FP64 numerical accuracy. Our solver is released as open-source software and is available at \url{https://github.com/Lhongpei/D-PDLP}.
\end{abstract}

\section{Introduction}
\label{sec:intro}
Large-scale linear programming (LP) lies at the core of numerous machine learning \citep{agarwal2018reductions,ben2013robust,taskar2005learning,joachims2009cutting}, operations research \citep{waissi1994network,gallego1997multiproduct,talluri2006theory}, and data-driven decision-making problems, including resource allocation \citep{hibiki2006multi,kelly1998rate}, market equilibrium computation \cite{eisenberg1959consensus,orlin2010improved}, and multi-stage stochastic programming \citep{hibiki2006multi,gangammanavar2021stochastic}. Modern applications,especially in machine learning, often involve millions or even billions of variables and constraints, rendering classical interior-point methods impractical due to their unfavorable memory footprint and limited parallel scalability.

First-order primal–dual methods, such as PDLP \citep{applegate2021practical}, cuPDLP.jl \citep{lu2023cupdlp}, cuPDLP-C \citep{lu2023cupdlpc}, cuPDLPx \citep{lu2024cupdlpx}, and HPR-LP \citep{chen2025hpr} have recently emerged as a promising alternative for solving large-scale LPs, owing to their low per-iteration complexity and amenability to parallelization. Among them, the PDLP family based on the primal–dual hybrid gradient (PDHG) framework has demonstrated strong empirical performance and robustness across a wide range of benchmark instances.

To date, the applicability of PDHG to huge-scale LP problems in multi-GPU environments remains unexplored. In this paper, we bridge this gap by proposing \textbf{D-PDLP}, the first \underline{D}istributed \underline{PDLP} framework. While distributed first-order methods historically suffer from prohibitive communication bottlenecks across loosely coupled nodes, D-PDLP is explicitly designed to exploit the massive intra-node bandwidth of modern dense GPU architectures (e.g., NVLink). Our approach spatially decomposes primal and dual updates across a grid of workers, executing partial Sparse Matrix-Vector Multiplications (SpMV) locally before synchronizing via high-speed \texttt{AllReduce} collectives. This tightly coupled design effectively mitigates communication bottlenecks, translating raw hardware bandwidth into algorithmic speedups.

Scaling PDHG effectively is non-trivial. Real-world LPs often exhibit irregular sparsity patterns that cause severe load imbalances under naive partitioning. To address this, we introduce a \textbf{block-wise random permutation} strategy combined with \textbf{nonzero-aware matrix partitioning}. Unlike naive fully random permutation, which destroys local data locality, our block-wise approach preserves dense micro-structures essential for efficient GPU computation while ensuring statistical load balance across the cluster.

Our main contributions are summarized as follows:
\begin{itemize}[leftmargin=0.5\leftmargin]
    \item We first propose a distributed PDHG-based solver for large-scale LPs, which is seamlessly extensible to quadratic programming (QP) and scales across multiple GPUs. The highly optimized C/CUDA implementation of D-PDLP is fully open-sourced.
    \item We propose a block-wise random shuffling and nonzero-aware partitioning strategy that effectively balances load under irregular sparsity while preserving local structure, which is of independent interest for large-scale sparse optimization.
    \item We establish a rigorous analytical cost model that mathematically characterizes the theoretical performance and scaling behavior of multi-GPU execution. This model is then utilized to automatically select the optimal GPU count prior to execution.
    \item Extensive experiments demonstrate that the proposed approach achieves strong multi-GPU scalability and high performance on both large-scale LP and diagonal QP benchmarks while maintaining full FP64 numerical accuracy. 
\end{itemize}

\subsection{Related Works} 
Traditional linear programming solvers, including the simplex method and interior-point methods (IPMs), are highly effective for small- to medium-scale problems. These methods deliver reliable, high-accuracy
solutions across various applications and are supported by commercial solvers such as MOSEK \citep{aps2019mosek}, GUROBI \citep{gurobi2021gurobi}, and COPT \citep{ge2022cardinal}. However, their applicability to large-scale instances is limited: the simplex method suffers from exponential worst-case complexity, while IPMs rely on sequential matrix factorizations that are difficult to scale and parallelize. As a result, first-order methods (FOMs) have emerged as a compelling alternative for solving large-scale LPs.

Among FOM-based approaches, ADMM-based solvers such as ABIP \citep{lin2021admm} and ABIP+ \citep{deng2022new} demonstrate the potential of first-order methods for large-scale LPs by requiring only a one-time matrix factorization. Building on this progress, PDHG-based solvers, most notably PDLP \citep{applegate2021practical} and its GPU implementation cuPDLP.jl \citep{lu2023cupdlp}, eliminate the need for matrix factorization altogether and represent the first methods capable of efficiently solving truly large-scale LPs. To further reduce runtime overhead, the C/CUDA-based implementation cuPDLP-C \citep{lu2023cupdlpc} was introduced. More recently, cuPDLPx \citep{lu2024cupdlpx} incorporates additional algorithmic heuristics and engineering optimizations to further improve convergence behavior and computational efficiency.

Beyond GPU-based first-order primal–dual methods developed specifically for LPs, a growing body of work has proposed GPU-accelerated solvers for a wide range of large-scale convex optimization problems. For quadratic programming, representative methods include rAPDHG \citep{lu2023practical}, PDHCG \citep{huang2025restarted}, and HPR-QP \citep{chen2025hpr}. For semidefinite programming, notable GPU-based solvers include LoRADS \citep{han2024low,han2024accelerating}, ALORA \citep{ding2025new}, and cuHALLaR \citep{aguirre2025cuhallar}. More generally, PDCS \citep{lin2025pdcs} addresses large-scale conic programming problems. GPU-enabled primal–dual methods have also been successfully applied to optimal transport \citep{lu2024pdot, zhang2025hot}, market equilibrium computation \citep{liu2025pdhcg}, and large-scale network flow problems \citep{zhang2025solving}. Collectively, these GPU-accelerated primal–dual solvers demonstrate substantial advantages over traditional approaches, particularly in terms of scalability and computational efficiency on modern parallel hardware.
\section{Preliminaries}
In this section, we briefly review the general formulation of linear programming in \Cref{sec:lp}, the primal–dual hybrid gradient  algorithm in \Cref{sec:pdhg}, and the cuPDLPx solver \citep{lu2024cupdlpx}, which serves as the foundation for our implementation, in \Cref{sec:cupdlpx}.
\subsection{Linear Programming}\label{sec:lp}
We consider the following linear programming problem 
\begin{equation}\label{prob:primal}
    \min_{x \in \mathcal{X}} \;c^\top x \quad \text{subject to } Ax \in \mathcal{S}
\end{equation}
where $A \in \mathbb{R}^{m \times n}$, $c \in \mathbb{R}^n$, the feasible set $\mathcal{X} := \{x \in \mathbb{R}^n : \ell_v \leq x \leq u_v\}$ with bounds $\ell_v \in (\mathbb{R} \cup \{-\infty\})^n$ and $u_v \in (\mathbb{R} \cup \{\infty\})^n$, and the constraint range $\mathcal{S} := \{s \in \mathbb{R}^m : \ell_c \leq s \leq u_c\}$
with $\ell_c \in (\mathbb{R} \cup \{-\infty\})^m$ and $u_c \in (\mathbb{R} \cup \{\infty\})^m$. One can readily derive the corresponding Lagrangian, which leads to an equivalent saddle-point formulation of the form
\begin{equation}\label{prob:saddle}
    \min_{x \in \mathcal{X}} \, \max_{y \in \mathcal{Y}} \; \mathcal{L}(x,y) \coloneqq c^Tx + y^T Ax -p(-y; \ell_c, u_c),
\end{equation}
where $p(\cdot \ ; l, u)$ denotes a piecewise linear function defined as
$$
p(y; l, u) = u^\top y^+ - l^\top y^-
$$
with $y^+ = \max\{y,0\}$ and $y^- = \max\{-y, 0\}$ being the positive and negative parts of $y$, respectively, and the dual feasible set $\mathcal{Y} \subseteq \mathbb{R}^m$ is Cartesian product whose $i$-th components are determined by the boundedness of the primal constraints:
\begin{equation*}
    \mathcal{Y}_i := 
    \begin{cases} 
    \{0\} & \text{if } (\ell_c)_i = -\infty, (u_c)_i = \infty, \\
    \mathbb{R}^- & \text{if } (\ell_c)_i = -\infty, (u_c)_i \in \mathbb{R}, \\
    \mathbb{R}^+ & \text{if } (\ell_c)_i \in \mathbb{R}, (u_c)_i = \infty, \\
    \mathbb{R} & \text{otherwise}.
    \end{cases}
\end{equation*}

\subsection{The PDHG Algorithm}\label{sec:pdhg}
The primal–dual hybrid gradient method can be applied to the saddle-point problem in \eqref{prob:saddle}, yielding the following update rules:
\begin{equation}
\begin{aligned}
    x^{k+1} &= \text{proj}_{\mathcal{X}} \left( x^k - \tau (c - A^\top y^k) \right) \\
    y^{k+1} &= y^k - \sigma z^{k+1} - \sigma \cdot \text{proj}_{-\mathcal{S}}\left(\sigma^{-1} y^k - A(2x^{k+1} - x^k) \right)
\end{aligned} 
\end{equation}
where $\tau$ and $\sigma$ denote the primal and dual step sizes, respectively. The step sizes can be reparameterized as  $\tau = \frac{\eta}{\omega}$ and $\sigma = \eta\omega$, where $\eta$ denotes the overall step size and $\omega$ is the primal weight that balances the primal and dual updates.

\subsection{Algorithmic Enhancements in cuPDLPx}\label{sec:cupdlpx}
Our distributed implementation builds upon the state-of-the-art GPU-based LP solver cuPDLPx \citep{lu2024cupdlpx}. To improve convergence behavior and numerical stability, cuPDLPx incorporates several algorithmic enhancements, including reflected Halpern iteration scheme, adaptive restarting strategy, and primal weight updates. These techniques further enhance the performance of PDLP. Detailed descriptions of these methods are provided in Appendix \ref{appd:ori_heuristic}.

\begin{figure}
\vspace{-1.5em}
    \centering
    \includegraphics[width=0.75\linewidth]{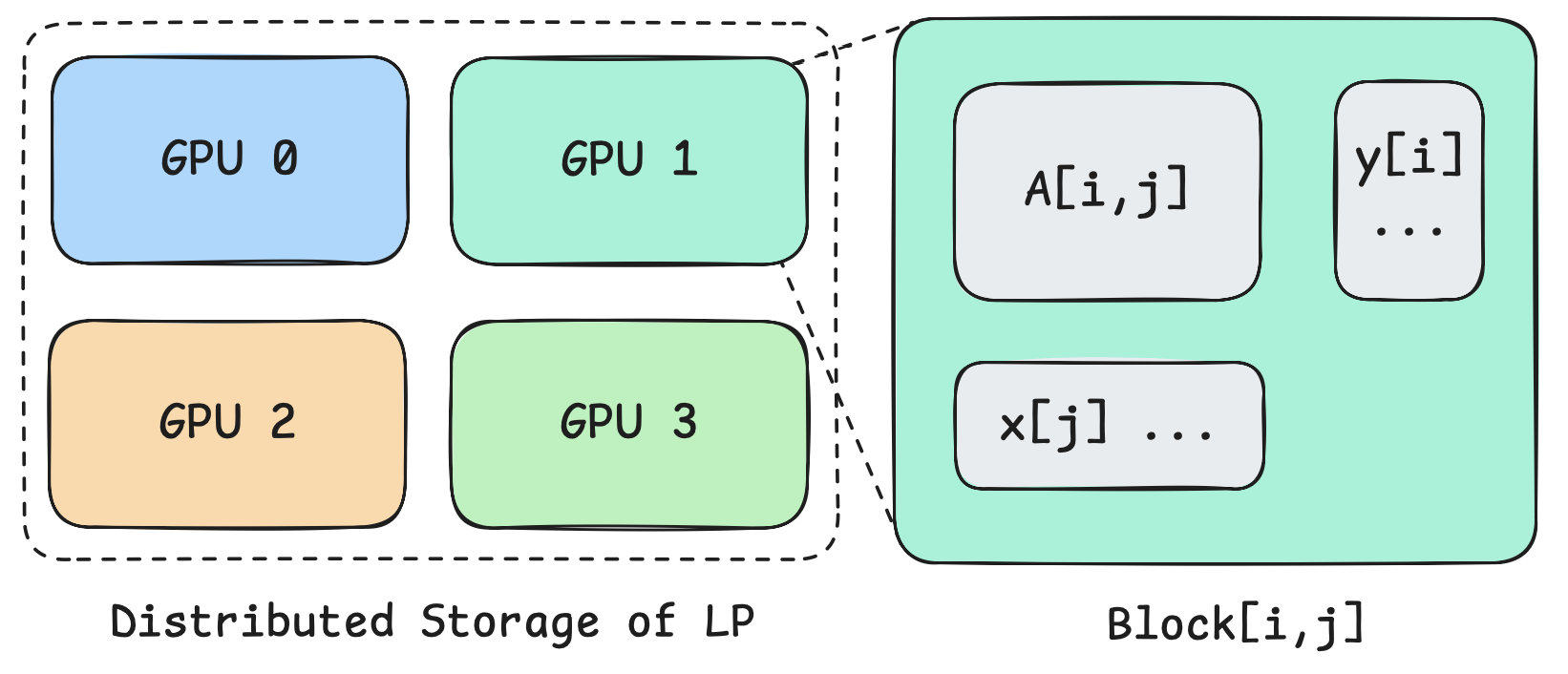}
    \caption{Distributed memory layout under 2D grid partitioning.
\textbf{Left:} A logical \(2 \times 2\) device mesh illustrating the global topology.
\textbf{Right:} Local data stored on the device at grid coordinate \((i,j)\). The constraint matrix is stored as block \(A_{[i,j]}\), primal variables \(x_{[j]}\), and dual variables \(y_{[i]}\). This layout enables local computation of \(A_{[i,j]} x_{[j]}\) and \(A_{[i,j]}^\top y_{[i]}\) prior to reduction.}
\label{fig:distributed_storage}
\vspace{-1.5em}
\end{figure}

\section{Distributed PDLP}
\label{sec:d-pdlp}
To scale PDHG beyond single-device computation and unlock the computational throughput of modern clusters, we adopt a 2D partitioning strategy for the constraint matrix $A$. This design enables the flexible distribution of variables and constraints across multiple GPUs, allowing the solver to effectively exploit massive parallelism to accelerate the solving phase.

\subsection{Problem Partition}
We adopt a variant of the named-axis notation to describe the data layout across hardware resources \citep{scaling-book}. Specifically, we assume a 2D grid of $N_{\mathrm{device}} = |R| \times |C|$ devices, where each axis is assigned a semantic role:
\begin{itemize}
\item \textbf{Row axis \((R)\)}: Corresponds to partitioning along the rows of $A$, with cardinality \(|R|\).
\item \textbf{Column axis \((C)\)}: Corresponds to partitioning along the columns of $A$, with cardinality \(|C|\).
\end{itemize}
A partition \(\mathcal{P}\) is defined by mapping the dimensions of a tensor to the corresponding physical grid axes. For the sparse constraint matrix \(A \in \mathbb{R}^{m \times n}\), applying the partition \(\mathcal{P}\) implies that the global row index set is divided into \(|R|\) disjoint subsets along the grid axis \(R\), while the global column index set is divided into \(|C|\) disjoint subsets along the grid axis \(C\). As a result, the matrix \(A\) is decomposed into a two-dimensional grid of local sparse blocks, as illustrated in Figure~\ref{fig:distributed_storage}:
\begin{equation} A \rightarrow
\begin{bmatrix}A_{[1,1]} & A_{[1,2]} & \dots & A_{[1, |C|]} \\A_{[2,1]} & A_{[2,2]} & \dots & A_{[2, |C|]} \\ \vdots & \vdots & \ddots & \vdots \\A_{[|R|, 1]} & A_{[|R|, 2]} & \dots & A_{[|R|, |C|]}
\end{bmatrix}
\end{equation}
where each sub-matrix $A_{[i,j]}$ contains the block of non-zero entries corresponding to the $i$-th row partition and the $j$-th column partition. Under this layout, each block $A_{[i,j]}$ is stored exclusively in the local memory of the device located at grid coordinate $(i, j)$.

\paragraph{Local Shapes and Sparsity.}
In dense tensor settings, as in LLM model parallelism such as Megatron-LM \citep{shoeybi2020megatronlmtrainingmultibillionparameter}, a partition \(\mathcal{P}\) typically assigns each device a uniform block of size \((m/|R|, \, n/|C|)\) \citep{scaling-book}. In contrast, our setting features irregular local shapes, where device $(i,j)$ stores submatrix $A_{[i,j]}$ of dimensions $(m_i, n_j)$, with $\sum_i m_i = m$ and $\sum_j n_j = n$. Moreover, the number of nonzero entries may vary significantly across devices. This irregularity motivates the load-balancing heuristics introduced in Section~\ref{heur_partition}. A deeper analysis of the differences between dense model parallelism and our framework is provided in Appendix~\ref{app:megatron}.

\paragraph{Vector Partitioning and Replication.}
To minimize communication during local computations, we align the distribution of all vectors with the two-dimensional partitioning of the constraint matrix. Primal-space vectors in \(\mathbb{R}^n\), including the decision variable \(x\) and objective coefficients \(c\), are partitioned along the column axis \(C\). The \(j\)-th block, denoted by \(x_{[j]}, c_{[j]} \in \mathbb{R}^{n_j}\), induces a local primal feasible set 
$$\mathcal{X}_{[j]} := \{z \in \mathbb{R}^{n_j} : (\ell_{v})_{[j]} \le z \le (u_{v})_{[j]}\}.$$
Each block is vertically replicated across all devices in column \(j\), ensuring that any device holding \(A_{[i,j]}\) can compute \(A_{[i,j]} x_{[j]}\) and \(\operatorname{proj}_{\mathcal{X}_{[j]}}\) without communication.

Similarly, dual-space vectors in \(\mathbb{R}^m\), including the dual variable \(y\), are partitioned along the row axis \(R\), with blocks \(y_{[i]} \in \mathbb{R}^{m_i}\) and corresponding local constraint sets $$
\mathcal{S}_{[i]} := \{b \in \mathbb{R}^{m_i} : (\ell_{c})_{[i]} \le b \le (u_{c})_{[i]}\}.
$$
These blocks are horizontally replicated across all devices in row \(i\), enabling local computation of \(A_{[i,j]}^\top y_{[i]}\) and dual projections \(\operatorname{proj}_{-\mathcal{S}_{[i]}}\). By orthogonally replicating primal and dual variables, communication is limited to reductions of partial matrix–vector products, maximizing the efficiency of local update kernels.

\subsection{Distributed PDHG Algorithm}
\label{sec:pdhg-iter}
The iterative updates necessitate global synchronization of partial matrix–vector products computed on local data shards. We implement this synchronization using the NCCL \texttt{AllReduce} collective. Formally, for a named grid axis \(X\), the operator \(\operatorname{AllReduce}_X(\cdot)\) aggregates partial vectors via summation across all devices along axis \(X\) and broadcasts the accumulated result back to them. The complete distributed procedure is outlined in Algorithm~\ref{alg:dist_pdlp}.

\subsubsection{Communication Operators}
We formalize the communication pattern by defining the operator \(\operatorname{AllReduce}_X\) acting on a distributed tensor \(T\), where \(T_{[i,j]}\) denotes the local block stored on the device at grid coordinate \((i,j)\). Depending on the aggregation axis, we define three variants:
\begin{align}
\operatorname{AllReduce}_R(T_{[i,j]}) = \sum_{k=1}^{|R|} T_{[k,j]}, \quad &\operatorname{AllReduce}_C(T_{[i,j]}) = \sum_{k=1}^{|C|} T_{[i,k]} \\
\operatorname{AllReduce}_G(T_{[i,j]}) &= \sum_{k=1}^{|R|} \sum_{l=1}^{|C|} T_{[k,l]}
\end{align}
Here, \(\operatorname{AllReduce}_R\) aggregates data across the row axis, \(\operatorname{AllReduce}_C\) across the column axis, and \(\operatorname{AllReduce}_G\) across the entire device grid. In our distributed PDHG implementation, \(\operatorname{AllReduce}_R\) is used to aggregate partial products in the computation of \(A^\top y\), while \(\operatorname{AllReduce}_C\) is used to sum partial results in the computation of \(Ax\).

\subsubsection{Distributed Primal–Dual Updates}
The primal variables \(x\) are partitioned along the column axis \(C\) into blocks \(\{x_{[j]}\}\), while the dual variables \(y\) are partitioned along the row axis \(R\) into blocks \(\{y_{[i]}\}\). 

For the primal update, the gradient term \(A^\top y^k\) is computed via partial products \(A_{[i,j]}^\top y_{[i]}^k\), followed by an \texttt{AllReduce} across the row axis (R):
\begin{equation*}
[A^\top y^k]_{[j]} = \operatorname{AllReduce}_R \bigl( A_{[i,j]}^\top y_{[i]}^k \bigr).
\end{equation*}
The primal block update is then given by
\begin{equation*}
x_{[j]}^{k+1} = \operatorname{proj}_{\mathcal{X}_{[j]}}
\left( x_{[j]}^k - \tau \bigl( c_{[j]} - [A^\top y^k]_{[j]} \bigr) \right).
\end{equation*}

For the dual update, let \(\bar{x} = 2x^{k+1} - x^k\). The partial products \(A_{[i,j]} \bar{x}_{[j]}\) are aggregated via an \texttt{AllReduce} across the column axis (C):
\begin{equation*}
z_{[i]} = \operatorname{AllReduce}_C \bigl( A_{[i,j]} \bar{x}_{[j]} \bigr),
\end{equation*}
followed by the local dual update
\begin{equation*}
y_{[i]}^{k+1} = y_{[i]}^k - \sigma z_{[i]}
- \sigma \cdot \operatorname{proj}_{-\mathcal{S}_{[i]}}
\bigl( \sigma^{-1} y_{[i]}^k - z_{[i]} \bigr).
\end{equation*}
This design ensures that the computationally intensive partial matrix-vector products \(A_{[i,j]} \bar{x}{[j]}\) and \(A_{[i,j]}^\top y_{[i]}\) are executed entirely within local memory without any inter-device communication. Synchronization is deferred until the aggregation phase, allowing all GPUs to saturate their compute throughput independently during the SpMV kernels.

Beyond the high-frequency core updates, infrequent operations such as KKT residual evaluations and heuristic enhancements (e.g., adaptive restart) require additional synchronization. Because they are executed at very low frequencies, their overall impact on the communication overhead is minimal (details in Appendix~\ref{appd:heuristic}).

\subsection{Per-iteration Time Complexity}
\label{subsec:per_iter_time}

Since auxiliary operations (e.g., KKT evaluations and enhancements) are executed infrequently, our theoretical cost model focuses strictly on the high-frequency core updates (primal, dual, and Halpern). Given the extremely low arithmetic intensity of the PDHG algorithm, the local execution is strictly memory-bandwidth bound (I/O bound) rather than compute-bound. This physical property significantly simplifies our analytical model: the execution time can be accurately quantified purely by the volume of memory traffic, bypassing complex instruction-level computations and negligible local kernel launch overheads.

To analytically model this I/O workload, we assume the non-zero elements ($\mathrm{nnz}$) are uniformly distributed across all local blocks. In practice, this uniformity is effectively guaranteed by our block-wise random permutation strategy, detailed in Section~\ref{heur_partition}.

Under these premises, let $N_{\mathrm{device}} = |R| \times |C|$ denote the total number of physical GPUs. The theoretical execution time of a single iteration, $T_{\mathrm{step}}$, is formulated as follows:

\begin{equation}
\label{eq:cost_model}
\begin{split}
T_{\mathrm{step}}(|R|, |C|) =& \frac{\mathcal{C}_{\mathrm{spmv}} \cdot \mathrm{nnz}}{B_{\mathrm{mem}} \cdot N_{\mathrm{device}}} + \frac{\mathcal{C}_{\mathrm{vec}}}{B_{\mathrm{mem}}} \left( \frac{n}{|C|} + \frac{m}{|R|} \right) \\
& + \frac{\mathcal{C}_{\mathrm{net}}}{B_{\mathrm{net}}} \left( \frac{n|R| + m|C| - (m+n)}{N_{\mathrm{device}}} \right) + n_{\mathrm{hop}} \beta_{\mathrm{sync}} \big( \mathbb{I}(|R| > 1) + \mathbb{I}(|C| > 1) \big)
\end{split}
\end{equation}

where $\mathrm{nnz}$ denotes the number of non-zero elements, and $m, n$ represent the dimensions of the constraint matrix. The hardware parameters $B_{\mathrm{mem}}$ and $B_{\mathrm{net}}$ represent the effective HBM and interconnect bandwidth, while $n_{\mathrm{hop}}$ and $\beta_{\mathrm{sync}}$ characterize the physical network hop count and the base collective synchronization latency, respectively. The coefficients $\mathcal{C}$ represent the bytes transferred per unit of data. Detailed theoretical derivations for each component, along with explicit definitions of all parameters, are deferred to Appendix~\ref{appd:theory}. 

Equation~\eqref{eq:cost_model} elucidates the physical basis of our distributed scalability: while the dominant SpMV memory I/O drops linearly ($\mathcal{O}(N_{\mathrm{device}}^{-1})$), the synchronization and vector I/O overheads under a balanced 2D grid scale only sub-linearly ($\mathcal{O}(N_{\mathrm{device}}^{-1/2})$). Beyond theoretical analysis, this model drives our \textbf{AUTO} scheduling heuristic. Although inherent GPU micro-architectural complexities make predicting exact absolute runtimes elusive, the model reliably captures relative scaling trends. By evaluating Eq.~\ref{eq:cost_model} prior to execution, AUTO adaptively selects the optimal GPU count. As confirmed in our experiments (Section~\ref{sec:experiments}), this analytical approach effectively improve D-PDLP performance across varying problem scales.

\section{Load-Balancing Strategies}
\label{heur_partition}
Achieving high parallel efficiency in distributed sparse optimization requires addressing three interconnected challenges: minimizing communication overhead relative to computation, mitigating load imbalance caused by irregular sparsity, and preserving memory locality for efficient GPU execution. We address these through a three-stage partition strategy: adaptive grid topology selection, block-wise random permutation, and nonzero-aware matrix partitioning. Crucially, these load-balancing strategies are computationally inexpensive and are performed only once during the preprocessing stage, leading to negligible overhead compared to the overall solving time. A detailed empirical analysis of this overhead is provided in Appendix~\ref{appd:overhead}.

\subsection{Adaptive Grid Topology Selection}
\label{sec:ada_grid}
To minimize the per-iteration communication overhead $\mathcal{O}(n/|C| + m/|R|)$, described in the Equation~\ref{eq:cost_model}, the $R \times C$ processor grid must dynamically adapt to the dimensions of the constraint matrix $A \in \mathbb{R}^{m \times n}$. Assuming uniform partitions via random permutation, we balance the primal and dual data transmissions by matching the grid's aspect ratio to the matrix's. Specifically, for a given device count $N_{\mathrm{device}}$, we configure the grid by selecting the factor pair $(|R|, |C|)$ that satisfies $|R| \times |C| = N_{\mathrm{device}}$ while minimizing the difference between $|R|/|C|$ and $m/n$.

\begin{figure*}
    \centering
    \includegraphics[width=0.90\linewidth]{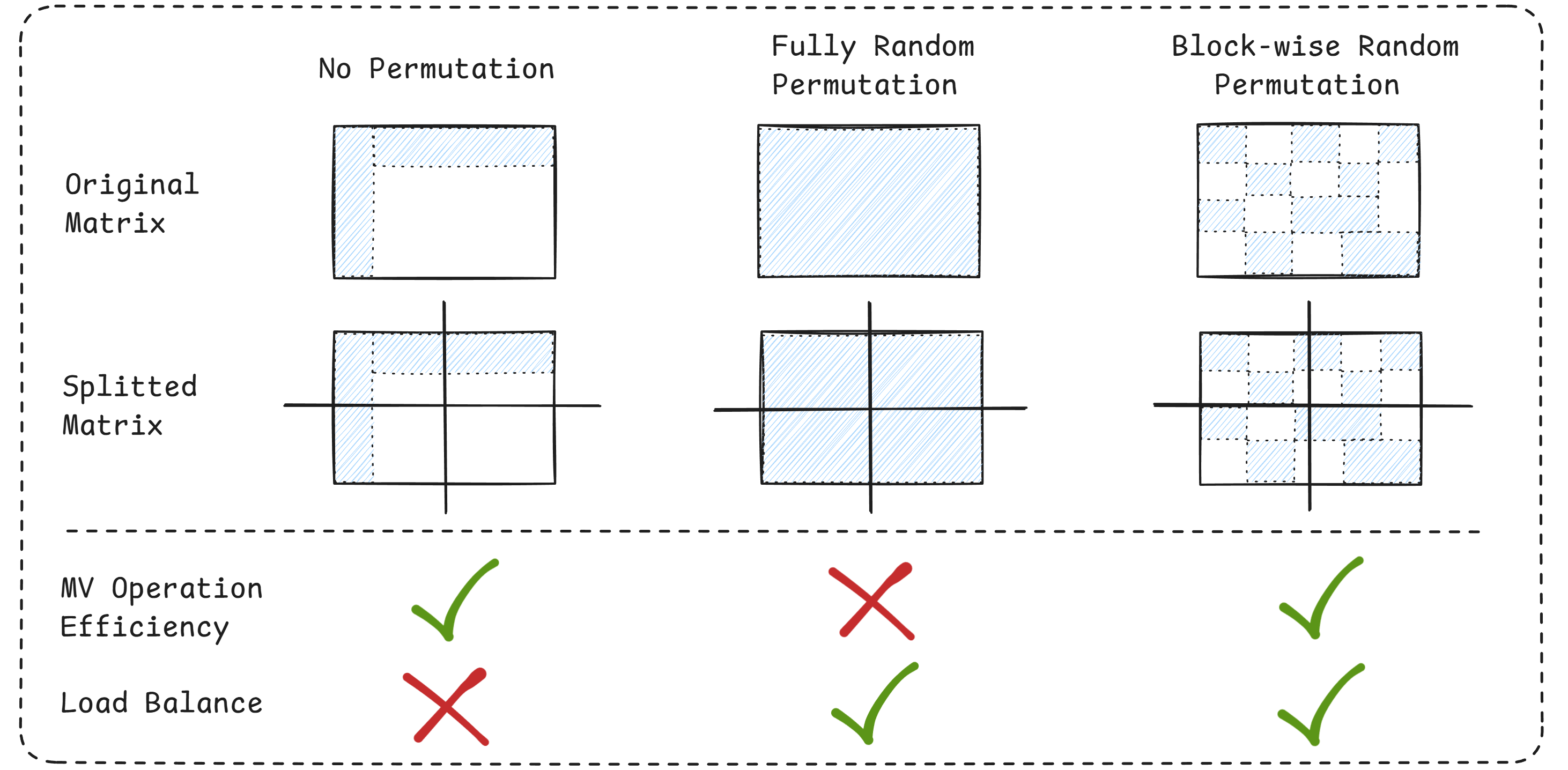}
    \caption{Comparison of matrix permutation strategies.
\textbf{Left:} Natural ordering preserves dense structures but causes severe load imbalance.
\textbf{Middle:} Full random permutation balances the workload but destroys local sparsity, degrading SpMV efficiency.
\textbf{Right:} Block-wise random permutation achieves both load balance and memory locality, enabling scalable GPU execution.}
\label{fig:permutation_strategy}
\vspace{-1.5em}
\end{figure*}

\subsection{Block-Wise Random Permutation}
\label{heur_permutation}
Irregular sparsity patterns in real-world LP instances, such as staircase structures or clustered nonzeros, can lead to severe computational load imbalance under a static 2D partition (see Figure~\ref{fig:permutation_strategy}, Left). While a \textit{Fully Random Permutation} of rows and columns can statistically homogenize the nonzero distribution (Figure~\ref{fig:permutation_strategy}, Middle), it introduces critical inefficiencies at the hardware level.
\paragraph{Inefficiency of Fully Random Permutation.}
Although fully random permutation achieves excellent global load balance, it destroys the inherent block-dense micro-structures common in industrial LPs. This fragmentation scatters nonzeros randomly in memory, leading to irregular access patterns that break memory coalescing. On GPUs, this results in significant thread divergence and low cache hit rates during SpMV operation, drastically degrading computational throughput.
\paragraph{Block-wise Random Permutation.}
To reconcile global load balancing with local SpMv efficiency, we introduce a \textit{Block-Wise Random Permutation} strategy (Figure~\ref{fig:permutation_strategy}, Right). Instead of permuting individual indices, we partition the global rows and columns into contiguous blocks of fixed size \(B\) (e.g., \(B=256\) to match CUDA thread blocks) and apply random permutations over these blocks. This approach disperses dense regions evenly across the \(|R| \times |C|\) device grid to ensure load balancing, while simultaneously preserving the dense micro-structures within each \(B \times B\) block to enable efficient SpMV execution on GPUs. Following empirical sensitivity analysis across diverse instances, we fixed $B=256$.

\subsection{Nonzero-Aware Matrix Partitioning}
Following permutation, we partition the matrix independently along each grid axis $X \in \{R, C\}$. We select boundaries $\{p_0, \ldots, p_{|X|}\}$ such that the nonzero count in each slice satisfies $\sum_{k=p_i}^{p_{i+1}-1} \text{nnz}_k \approx \text{nnz}(\tilde{A}) / |X|$. While this axis-wise heuristic alone cannot guarantee strict 2D load balance, its synergy with block-wise random permutation effectively ensures a uniform workload distribution.

\section{Numerical Experiments}
In this section, we present a comprehensive evaluation of the proposed D-PDLP framework, assessing its scalability, computational efficiency, and ability to harness multi-GPU utility to accelerate PDLP.
\label{sec:experiments}
\subsection{Experimental Setup}

We evaluate our distributed solver on a single dense compute node equipped with eight NVIDIA H100 (80GB) GPUs fully interconnected via NVLink. This tightly coupled topology is essential for minimizing the communication overhead of D-PDLP's dense \texttt{AllReduce} operations. For maximum performance, the solver is implemented in C/CUDA extending \texttt{cuPDLPx} \citep{lu2024cupdlpx}, alongside a Julia prototype for architecture validation. Comprehensive hardware specifications and implementation details are provided in Appendix~\ref{appd:exp_setup}.

\textbf{Evaluation Metric.}\quad We aggregate performance across instances using the standard Shifted Geometric Mean (SGM) with a shift of $k=10$ seconds~\citep{mittelmann2020benchmarks}, defined as $\text{SGM10} = \exp \big( \frac{1}{N} \sum_{i=1}^N \ln(t_i + 10) \big) - 10$, where $t_i$ is the solving time for the $i$-th instance.

\textbf{Automatic GPU Count (AUTO):}\quad To optimize the performance of D-PDLP across varying problem scales, we implement an automatic scheduler, AUTO, which selects the optimal GPU count ($N \in \{1, 2, 4, 8\}$) prior to execution. AUTO is guided by the theoretical cost model defined in Eq.~\ref{eq:cost_model}, utilizing hardware-specific parameters ($B_{\mathrm{mem}}=3350$~GB/s, $B_{\mathrm{net}}=400$~GB/s, and $\beta_{\mathrm{sync}}=0.01$~ms) to estimate the per-step cost for a given instance ($m, n, \mathrm{nnz}$). Detailed parameter estimations and examples are provided in Appendix~\ref{appd:theory}.

\begin{table}[htbp]
    \vspace{-1em}
    \centering
    \caption{Performance on standard benchmarks. Time is reported as SGM10. `Solved' indicates instances solved within the time limit (3600s for Mittelmann, 1200s for MIPLIB). Best results are \textbf{bold}, second-best are \underline{underlined}.}
    \label{tab:sgm_results_merged}
    \resizebox{0.6\linewidth}{!}{%
    \begin{tabular}{llcccccc}
        \toprule
         & & \multicolumn{5}{c}{Number of GPUs} & \\
        \cmidrule(l){3-7}
        Dataset & Metric & 1 & 2 & 4 & 8 & AUTO \\
        \midrule
        \multirow{2}{*}{Mittelmann} & Solved & \textbf{46} & \textbf{46} & \textbf{46} & \textbf{46} & \textbf{46} \\
         & Time (s) & 21.26 & 19.00 & \underline{18.58} & 19.24 & \textbf{18.38} \\
        \midrule
        \multirow{2}{*}{\shortstack[l]{MIPLIB Small\\(100K - 1M)}} & Solved & \textbf{264} & \textbf{264} & \textbf{264} & \underline{263} & \textbf{264} \\
         & Time (s) & 5.82 & \textbf{5.77} & 6.12 & 6.79 & \underline{5.81} \\
        \midrule
        \multirow{2}{*}{\shortstack[l]{MIPLIB Med.\\(1M - 10M)}} & Solved & \textbf{92} & \textbf{92} & \textbf{92} & \textbf{92} & \textbf{92} \\
         & Time (s) & 7.21 & \textbf{6.68} & \underline{6.70} & 7.11 & 6.76 \\
        \midrule
        \multirow{2}{*}{\shortstack[l]{MIPLIB Large\\($>10$M)}} & Solved & 13 & \underline{14} & \textbf{20} & \textbf{20} & \textbf{20} \\
         & Time (s) & 90.34 & 75.30 & \underline{62.07} & \textbf{56.81} & \textbf{56.81} \\
        \midrule
        \multirow{2}{*}{\shortstack[l]{MIPLIB Overall}}& Solved & 369 & 370 & \textbf{376} & \underline{375} & \textbf{376} \\
         & Time (s) & 7.80 & \underline{7.47} & 7.60 & 8.14 & \textbf{7.31} \\
        \bottomrule
    \end{tabular}%
    }
\end{table}

\paragraph{Benchmark Datasets}
To rigorously evaluate performance and scalability, we utilize two distinct categories of benchmark datasets, ranging from standard academic problems to massive industrial-scale instances.

\begin{itemize}[leftmargin=0.2\leftmargin]
    \item \textbf{Standard Benchmarks.}\quad For initial validation and strong scaling analysis, we employ the MIPLIB 2017 \citep{gleixner2021miplib} relaxation set, comprising 383 diverse mixed-integer programming instances solved as linear programs, alongside Mittelmann’s LP benchmark \citep{mittelmann2005decision}, which contains 49 challenging public instances. We restrict our focus to feasible and bounded instances, setting the termination tolerance to $\epsilon = 10^{-4}$.
    \item \textbf{Huge-Scale Instances.}\quad Moving beyond standard benchmarks, we probe the architecture using a curated, huge-scale dataset. This collection includes the \texttt{zib03} instance from~\citep{kochProgressMathematicalProgramming2022}, large-scale PageRank formulations \citep{nesterov2014subgradient} ($10^7$ nodes), and real-world network datasets (e.g., LiveJournal) \citep{leskovec2014snap}. We also incorporate extremely large multicommodity flow (\texttt{mcf}) and synthetic design matching (\texttt{sdm}) instances. Furthermore, the suite includes Quadratic Assignment Problem (QAP) relaxations based on Adams-Johnson linearization \citep{adams1994improved} (such as \texttt{tai50b}, \texttt{lipa50b}) and Unit Commitment problems (\texttt{ds1}, \texttt{ds2}). For these instances, we adopt a tighter accuracy of $\epsilon = 10^{-6}$ (except \texttt{mcf} at $10^{-4}$).
\end{itemize}

\subsection{Performance on Benchmark Datasets}
\label{subsec:performance}

Tables~\ref{tab:sgm_results_merged} and~\ref{tab:performance} demonstrate the scalability of our distributed, grid-partitioned solver. For the smallest problem set, specifically the \textit{MIPLIB Small} subset with fewer than one million nonzeros, performance peaks early at just 2 GPUs (5.77s) and noticeably degrades when scaling to 8 GPUs (6.79s). This behavior is expected: for small problems, local GPU resources are easily saturated. Consequently, the marginal computational gains from wider parallelism cannot offset the escalating communication overhead introduced by dense \texttt{AllReduce} operations across the network.

This scaling saturation on trivial instances is acceptable, as first-order solvers are inherently designed to target large- and huge-scale problems. This design goal is confirmed by the strong scaling observed on larger instances. In the \textit{MIPLIB Large} category, the distributed implementation reduces SGM10 from 90.34s on a single GPU to 56.81s on 8 GPUs, demonstrating substantial performance gains from multi-GPU execution. Crucially, our AUTO scheduler successfully navigates these diverse architectural trade-offs without empirical profiling. By circumventing the latency wall on small instances while fully exploiting parallelism for massive problems, AUTO secures the fastest SGM10 times across both benchmarks, achieving 18.38s on Mittelmann and an overall best of 7.31s on MIPLIB, , demonstrating the practical effectiveness of our theoretical cost model.

The advantage is decisive on huge-scale instances (Table~\ref{tab:performance}). We observe \textbf{near-linear acceleration} on computationally intensive problems: \texttt{mcf\_2500} achieves speedups of $3.1\times$ (4 GPUs) and $6.0\times$ (8 GPUs), while \texttt{sdm\_50k} attains $6.2\times$ on 8 GPUs. Even the highly structured \texttt{zib03} sees a $3.3\times$ boost. Guided by the theoretical cost model, our AUTO scheduler correctly predicts the optimal 8-GPU configuration across all such instances, which perfectly aligns with the physical intuition that extreme-scale problems are heavily compute-bound. These results confirm that the grid-partitioned architecture effectively translates additional compute resources into proportional reductions in solving time. An end-to-end comparison with representative external solvers is provided in Appendix~\ref{app:external}.

\textbf{Improved Implementation.}\quad To further minimize runtime overhead, our solver is ultimately implemented in C/CUDA extending \texttt{cuPDLPx}~\citep{lu2024cupdlpx}. This optimized implementation reduces the aggregate SGM10 solving time by roughly 15\% to 20\% compared to our Julia prototype; a comprehensive instance-by-instance comparison is deferred to Appendix~\ref{appd:c_optimization}.

\subsection{Ablation Study on Load-Balancing Strategies}
To isolate the performance gains introduced by our load-balancing techniques, we evaluated the individual and combined effects of matrix permutation and partitioning strategies. As detailed in Appendix~\ref{appd:ablation_permute}, the baseline approach (\textit{No Permutation} combined with \textit{No NNZ} uniform partitioning) suffers from severe load imbalance, resulting in an aggregate SGM10 of 220.08s. In contrast, combining our proposed \textit{Block Random} permutation with \textit{Yes NNZ} nonzero-aware partitioning optimally disperses dense regions while preserving local block structures for efficient SpMV execution. This synergistic configuration reduces the SGM10 to 75.85s, delivering the most robust performance across the benchmark suite. Accordingly, all subsequent experiments adopt this configuration.

\begin{table*}[tbp]
\vspace{-0.5em}
    \centering
    \caption{Performance comparison: Single GPU vs. Distributed PDLP on massive-scale instances. Values represent wall-clock time in seconds. Best results are marked in \textbf{bold}, and second-best are \underline{underlined}. Notably, our AUTO scheduler correctly predicts the optimal 8-GPU configuration for all instances listed.}
    \label{tab:performance}
    \resizebox{0.9\textwidth}{!}{%
    \begin{tabular}{llrrr r rr}
        \toprule
        & & & & &  \multicolumn{3}{c}{Number of GPUs} \\
        \cmidrule(l){6-8}
        Source & Instance & $m$ & $n$ & nnz & 1& 4& 8 (AUTO)\\
        \midrule
        Koch & zib03 & 19,731,970 & 29,128,799 & 104,422,573 & 812 & \underline{336} & \textbf{245} \\
        \midrule
        \multirow{3}{*}{Pagerank} & rand-10m-nodes & 10,000,001 & 10,000,000 & 79,999,982 & 10 & \underline{5} & \textbf{4} \\
         & com-livejournal & 3,997,963 & 3,997,962 & 77,358,302 & 3 & \underline{2} & \textbf{1} \\
         & soc-livejournal1 & 4,847,572 & 4,847,571 & 78,170,533 & 2& \underline{1}& \textbf{1} \\
        \midrule
        \multirow{3}{*}{\shortstack{multicommodity-\\flow}} & mcf\_2500\_100\_500 & 1,512,600 & 126,250,100 & 253,750,100 & 2943 & \underline{935} & \textbf{504} \\
         & mcf\_5000\_50\_500 & 2,775,050 & 126,250,050 & 253,750,050 & 9114 & \underline{2994} & \textbf{1677} \\
         & mcf\_5000\_100\_250 & 1,775,100 & 127,500,100 & 257,500,100 & 9173 & \underline{3159} & \textbf{1732} \\
        \midrule
        design-matching & sdm\_50k\_500k\_15\_10 & 5,500,135 & 10,000,000 & 690,000,000 & 377 & \underline{92} & \textbf{61} \\
        \midrule
        \multirow{5}{*}{QAP} & wil50 & \multirow{5}{*}{3,437,600} & \multirow{5}{*}{6,252,500} & \multirow{5}{*}{19,125,000} & 43 & \underline{31} & \textbf{25} \\
         & lipa50a & & & & 27 & \underline{20} & \textbf{17}\\
         & lipa50b & & & & 28 & \underline{21} & \textbf{17} \\
         & tai50a & & & & 30 & \underline{21} & \textbf{18} \\
         & tai50b & & & & 82 & \underline{52} & \textbf{43} \\
        \midrule
        \multirow{2}{*}{Unit Com.} & ds1 & \multirow{2}{*}{641,037} & \multirow{2}{*}{659,145} & \multirow{2}{*}{21,577,566} & 279 & \underline{143} & \textbf{122} \\
         & ds2 & & & & 1288 & \underline{691} & \textbf{591} \\
         \midrule
        \multicolumn{5}{l}{\textbf{SGM10}} & 186.6 & \underline{101.2} & \textbf{77.7} \\
        \bottomrule
    \end{tabular}%
    }
    \vspace{-1.3em}
\end{table*}

\section{Extensions to Diagonal QPs: D-PDHCG}
\label{sec:extension_qp}

While D-PDLP is primarily designed for linear programming, the framework naturally extends to convex quadratic programming (QP) problems featuring a diagonal quadratic objective matrix $Q$. In this scenario, the primal update step is modified as follows:
\begin{equation}
x^{k+1} = \arg\min_{x \in \mathcal{X}} \left\{ \frac{1}{2}x^\top Q x + c^\top x + x^\top A^\top y^k + \frac{1}{2\tau_k} \|x - x^k\|_2^2 \right\}.
\label{eq:primal_update}
\end{equation}

Crucially, because $Q$ is diagonal, this subproblem remains fully separable across the dimensions of $x$ and admits a highly efficient closed-form solution. Under our distributed memory layout, the diagonal entries of $Q$ are partitioned along the column axis in exact alignment with the primal variable $x_{[j]}$. As a result, the primal update can be executed within local device memory, introducing zero additional communication overhead compared to the LP version.

To demonstrate the efficacy of this extension, Table~\ref{tab:real_lasso_full} presents an end-to-end performance comparison on large-scale LASSO instances against several state-of-the-art QP solvers, including HPR-QP, PDHCG, PDQP, SCS (GPU), OSQP, and COPT. As the results indicate, many traditional and first-order solvers struggle with memory constraints or exceed the time limit (denoted by ``t'', $>7200$s) on the largest instances, particularly those in the \texttt{avazu} and \texttt{kdd} families. In contrast, D-PDHCG successfully leverages multi-GPU parallelism to achieve the fastest solving times across all evaluated instances. This performance gap underscores the practical advantage of our zero-communication diagonal objective partitioning, which enables near-linear speedups on massive problems that are otherwise intractable for standard solvers. Detailed dataset descriptions, solver references, and experimental settings are provided in Appendix~\ref{appd:qp_dataset}.

\begin{table}[htbp]
  \centering
  \renewcommand{\arraystretch}{1.5}
  \caption{Solving time (s) over large-scale LASSO problems. The best results are marked in \textbf{bold}, and the second best are \underline{underlined}. ``t'' denotes time limit exceeded ($>7200$s).}
  \medskip
  \resizebox{\textwidth}{!}{
  \begin{threeparttable}
    \begin{tabular}{c|c|ccc|cccccc}
    \toprule
    Problem & \textbf{PDHCG-II}& \multicolumn{3}{c|}{\textbf{D-PDHCG (Ours)}} & HPR-QP & PDHCG & PDQP & SCS(GPU) & OSQP & COPT \\
            &                     & 2 GPUs & 4 GPUs & 8 GPUs &        &        &        &          &        &        \\
    \midrule
    SLS            & 1.96   & 1.10   & \underline{0.56}   & \textbf{0.50}   & 2.09   & 3.35   & 7.30   & 345.09 & 80.32  & 88.21  \\
    rcv1\_test     & 2.37   & 1.48   & \underline{1.06}   & \textbf{0.81}   & 6.21   & 7.12   & 19.54  & t      & t      & t      \\
    avazu-site.tr  & 1337.05& 953.02 & \underline{547.28} & \textbf{376.90} & 4911.41& 1377.54& 5124.82& t      & t      & t      \\
    avazu-app      & 217.70 & 116.92 & \underline{65.86}  & \textbf{46.30}  & 753.65 & 1429.55& 5557.97& t      & t      & t      \\
    avazu-site     & 1642.69& 993.62 & \underline{765.46} & \textbf{503.99} & 3213.38& 4224.95& t      & t      & t      & t      \\
    kddb2010\_test & 11.59  & 11.18  & \underline{8.79}   & \textbf{7.07}   & 26.87  & 11.10  & 46.49  & 490.66 & 255.81 & 69.57  \\
    kdda2010\_test & 9.90   & 8.50   & \underline{6.77}   & \textbf{5.61}   & 29.36  & 61.00  & 148.01 & t      & t      & t      \\
    kddb2010\_train& 703.36 & 466.28 & \underline{317.69} & \textbf{233.44} & 1971.24& 387.00 & 1715.50& t      & t      & t      \\
    kdda2010\_train& 341.16 & 237.19 & \underline{167.62} & \textbf{124.89} & 842.02 & 2705.94& t      & t      & t      & t      \\
    \bottomrule
    \end{tabular}
  \end{threeparttable}
  }
  \label{tab:real_lasso_full}
\end{table}
\section{Conclusion, Limitation, and Future Work}
\label{sec:conclusion}




We introduced \textbf{D-PDLP}, a distributed PDHG framework for solving large-scale LPs on multi-GPU systems. By combining two-dimensional grid partitioning with block-wise random permutation and nonzero-aware partitioning, D-PDLP balances irregular sparse workloads while preserving local structures for efficient SpMV. Experiments on large-scale benchmark and industrial instances demonstrate strong scalability, achieving up to $6\times$ speedups on 8 GPUs.

The tightly coupled design of D-PDLP is best suited to high-bandwidth intra-node architectures, where dense AllReduce communication can be efficiently supported. Although scaling beyond a single chassis remains limited by inter-node latency and bandwidth, emerging dense GPU systems such as NVIDIA GB200 NVL72 further expand the practical scope of this approach.

Additionally, while our theoretical cost model effectively guides resource allocation, predicting exact absolute execution times remains challenging due to micro-architectural complexities. Future work will focus on refining this analytical model to bridge the gap between theoretical predictions and real-world execution. Furthermore, we aim to extend the D-PDLP framework to a broader class of convex optimization problems, including Second-Order Cone Programming (SOCP) and Semidefinite Programming (SDP).

\newpage

\bibliographystyle{plainnat}
\bibliography{main,ref}

\newpage
\appendix
\onecolumn


\appendix

\section{Algorithm of Distributed PDLP}
We summarize the full initialization process and the distributed iteration logic in Algorithm~\ref{alg:dist_pdlp}.
\begin{algorithm*}[b]
  \caption{D-PDLP}
  \label{alg:dist_pdlp}
  \begin{algorithmic}
    \STATE {\bfseries Input:} Problem instance $(A, b, c)$, device count $N_{\mathrm{device}} \geq 2$, permutation $\Pi$, partition $\mathcal{P}$, initial point $(x^0, y^0)$, step sizes $\tau, \sigma$, and accuracy $\epsilon$.
    \STATE {\bfseries Initialization:} Determine grid topology $|R| \times |C| = N_{\mathrm{device}}$. Apply permutation $\Pi$ and partition $\mathcal{P}$ to get $$A \to \left\{A_{[i, j ]} \right\}, \; b \to \{b_{[i]}\}, \; c \to \{c_{[j]}\}, \; y^0 \to \{y^0_{[i]}\}, \; x^0 \to \{ x^0_{[j]}\}.$$ Store local blocks $(A_{[i,j]}, b_{[i]}, c_{[j]}, x_{[j]}^0, y_{[i]}^0)$ in each $(i,j)$-th GPU.
    \FOR{$k=0, 1, \cdots $, on each $(i,j)$-th GPU}
    
    \STATE \textbf{1. Distributed Primal Update}
    \STATE Compute local partial gradient: $g_{[i,j]} = A_{[i,j]}^\top y_{[i]}^k$
    \STATE Communication: Aggregate $g_{[i,j]}$ across row axis $R$: $[A^\top y^k]_{[j]} = \operatorname{AllReduce}_R \left( g_{[i,j]} \right)$
    \STATE Update local primal variable: $x_{[j]}^{k+1} = \operatorname{proj}_{\mathcal{X}_{[j]}} \left( x_{[j]}^k - \tau (c_{[j]} - [A^\top y^k]_{[j]}) \right)$
    \STATE \textbf{2. Distributed Dual Update}
    \STATE Compute local partial constraint: $v_{[i,j]} = A_{[i,j]} (2x_{[j]}^{k+1} - x_{[j]}^k)$
    \STATE Communication: Aggregate $v_{[i,j]}$ across column axis $C$: $z_{[i]} = \operatorname{AllReduce}_C \left( v_{[i,j]} \right)$
    \STATE Update local dual variable: $y_{[i]}^{k+1} = y_{[i]}^k - \sigma z_{[i]} - \sigma \cdot \operatorname{proj}_{-S_{[i]}} \left( \sigma^{-1} y_{[i]}^k - z_{[i]} \right)$
    \STATE \textbf{3. Distributed Halpern Iteration Update}
    \STATE \textbf{4. Termination criterion: } Stop if the KKT error is less than $\epsilon$.
    \ENDFOR
  \end{algorithmic}
\end{algorithm*}

The Halpern iteration scheme computes a linear combination of the current PDHG output, the previous iterate $z^k$, and the anchor point $z^0$. Since $z^0$ follows the identical partitioning scheme as the current iterate, this update is executed locally on each device without communication. Let $z_{[i,j]}^k = (x_{[j]}^k, y_{[i]}^k)$ denote the local portion of the primal-dual pair on device $(i, j)$. The update is:
\begin{equation}
    z_{[i,j]}^{k+1} = \left( (1+ \gamma )\frac{k+1}{k+2} \operatorname{PDHG}(z_{[i,j]}^k) - \gamma z_{[i,j]}^k \right) + \frac{1}{k+2} z_{[i,j]}^0
\end{equation}
This strict locality property ensures that the acceleration step introduces zero inter-device communication overhead.

\section{Algorithm Enhancement and Evaluation in cuPDLPx} \label{appd:ori_heuristic}
The convergence of PDHG is analyzed using the canonical norm $\|\cdot\|_P$, where the metric matrix $P$ is defined as:
\begin{equation}
    P := P_{\eta,\omega} = \begin{bmatrix} \frac{\omega}{\eta}I & A^\top \\ A & \frac{1}{\eta\omega}I \end{bmatrix}
\end{equation}
The iterate update is denoted as $z^{k+1} = \text{PDHG}(z^k)$, where $z^k = [x^k; y^k]^\top$.
\subsection{Reflected Halpern Iteration Scheme}
Rather than using vanilla PDHG iterates, \cite{lu2024restarted} adopts the restarted Halpern PDHG (rHPDHG) scheme. This scheme interpolates between the standard PDHG iterate and an initial anchor point $z^0$:
\begin{equation}
    z^{k+1} = \left( (1+ \gamma )\frac{k+1}{k+2} \text{PDHG}(z^k) - \gamma z^k\right) + \frac{1}{k+2} z^0
\end{equation}
This allows the algorithm to take more aggressive steps, resulting in stronger empirical performance \cite{lu2024cupdlpx}.

\subsection{Adaptive Restarting Strategy}
The adaptive restart strategy is designed to ensure the anchor point remains relevant to the local geometry. A restart is triggered based on a fixed-point error metric $r(z) = \|z - \text{PDHG}(z)\|_P$. The criteria for triggering a restart include:
\begin{itemize}
    \item \textbf{Sufficient Decay:} $r(z^{n,k}) \le \beta_{\text{sufficient}} r(z^{n,0})$.
    \item \textbf{Necessary Decay:} $r(z^{n,k}) \le \beta_{\text{necessary}} r(z^{n,0})$ and $r(z^{n,k}) > r(z^{n,k-1})$.
    \item \textbf{Artificial Restart:} $k \ge \beta_{\text{artificial}} T$, where $T$ is the total iterations.
\end{itemize}

\subsection{PID-Controlled Primal Weight Update}
To balance progress between the primal and dual spaces, the primal weight $\omega$ is dynamically adjusted using a Proportional-Integral-Derivative (PID) controller at each restart~\cite{lu2024restarted}. Logarithmic error $e^n$ is defined as the distance between the primal and dual to the anchor:
\begin{equation}
    e^n = \log \left( \frac{\sqrt{\omega^n} \|x^{n,t} - x^{n,0}\|_2}{\frac{1}{\sqrt{\omega^n}} \|y^{n,t} - y^{n,0}\|_2} \right)
\end{equation}
The weight for the next epoch is updated as:
\begin{equation}
    \log \omega^{n+1} = \log \omega^n - \left[ K_P \cdot e^n + K_I \sum_{i=1}^n e^i + K_D (e^n - e^{n-1}) \right]
\end{equation}
where $K_P$, $K_I$, and $K_D$ are the controller coefficients.

\subsection{KKT Evaluation}
The dual form of the problem is given by 
\begin{equation}
    \max_{y \in \mathcal{Y}, s \in \mathcal{R}} -p(-y; \ell_c, u_c) - p(-r; \ell_v, u_v) \quad \text{subject to } c - A^\top y = r,
\end{equation}
where $\mathcal{R} \subseteq \mathbb{R}^n$ is the Cartesian product whose $i$-th components are determined by the boundedness of the primal constraints:
\begin{equation}
    \mathcal{R}_i := 
    \begin{cases} 
    \{0\} & \text{if } (\ell_v)_i = -\infty, (u_v)_i = \infty, \\
    \mathbb{R}^- & \text{if } (\ell_v)_i = -\infty, (u_v)_i \in \mathbb{R}, \\
    \mathbb{R}^+ & \text{if } (\ell_v)_i \in \mathbb{R}, (u_v)_i = \infty, \\
    \mathbb{R} & \text{otherwise.}
    \end{cases}
\end{equation}
The solver terminates when the iterates satisfy a specified tolerance $\epsilon$ across three relative KKT error metrics: primal residual ($r_{\text{primal}}$), dual residual ($r_{\text{dual}}$), and the duality gap ($r_{\text{gap}}$).
\begin{align*}
    r_{\text{primal}} &= \frac{ \|Ax - \text{proj}_{\mathcal{S}}(Ax)\|_2}{1 + \|[\ell_c, u_c]\|_2}, \\[6pt]
    r_{\text{dual}} &= \frac{\left\| \text{proj}_{ \mathcal{X}}(x - \tau(c - A^\top y) ) - x \right\|_2 / \tau}{1 + \|c\|_2}, \\[6pt]
    r_{\text{gap}} &= \frac{|c^\top x + p(-y; \ell_c, u_c) - s^\top x|}{1 + \max\{|c^\top x|, |-p(-y; \ell_c, u_c) + s^\top x|\}}.
\end{align*}
where \(s = (\text{proj}_{ \mathcal{X}}(x - \tau(c - A^T y) ) - (x - \tau(c - A^T y))) / \tau \). The overall residual is defined as $\max\{r_{\text{primal}}, r_{\text{dual}}, r_{\text{gap}}\}$. 

\section{Differences from Dense Model-Parallel Systems}
\label{app:megatron}

Dense model-parallel training systems, exemplified by Megatron-LM \citep{shoeybi2020megatronlmtrainingmultibillionparameter}, partition large weight tensors across devices to scale deep-learning workloads. Their 2D grid layouts share a superficial geometric similarity with the partitioning scheme used in D-PDLP, yet the underlying computational settings and systems constraints differ fundamentally. Below we summarize the key distinctions.

\paragraph{Dense vs. Sparse Workloads.}
In dense model parallelism, each local block has a perfectly regular shape (e.g., $m/|R| \times n/|C|$) and the workload is statically predictable. In contrast, LP constraint matrices are highly sparse and structurally irregular. A naive row--column split can produce orders-of-magnitude imbalance in nonzero counts across devices, which directly translates to stragglers and under-utilized GPUs. Our block-wise random permutation together with nonzero-aware partitioning, illustrated in Section~\ref{heur_partition} is introduced precisely to address this irregularity.

\paragraph{Solver-Specific Design for PDLP.}
D-PDLP is designed around the PDLP iteration structure rather than being a generic sparse-matrix partitioning scheme. In each iteration, the two SpMV operations ($A\bar{x}$ and $A^{\top}y$) can be performed locally on each device once the corresponding primal and dual slices are available, while all element-wise operations (primal/dual updates, projections, and Halpern steps) require no inter-device communication. To support this execution pattern, the 2D layout is co-designed with the algorithm: (i) primal slices $x_{[j]}$ are replicated across row communicators so that each device in column $j$ can compute its local $A_{[i,j]}\bar{x}_{[j]}$ independently; (ii) dual slices $y_{[i]}$ are replicated across column communicators so that each device in row $i$ can compute its local $A_{[i,j]}^{\top}y_{[i]}$ independently; and (iii) two orthogonal AllReduce phases synchronize the partial results---one along columns for $A\bar{x}$ and one along rows for $A^{\top}y$. Consequently, the data layout, variable replication pattern, and communication schedule are all dictated by the PDLP algorithmic structure.

\section{Distributed Infrequent Operations and Communication Analysis}
\label{appd:heuristic}

In our distributed PDHG implementation, inter-device communication plays a critical role in determining overall system scalability. As discussed in the main text, the high-frequency core loop (primal, dual, and Halpern updates) strictly dominates the per-iteration runtime. 

However, a complete distributed solver requires additional operations, such as convergence monitoring and heuristic enhancements (e.g., adaptive restarts and PID weight updates adopted from cuPDLPx). These operations require exact matrix-vector products rather than intermediate differential quantities, necessitating heavy global synchronizations on the $R \times C$ processor grid. Because they are executed at remarkably low frequencies (e.g., strictly every $K$ iterations, where $K$ is typically large), their amortized communication and computation overheads are negligible. Table~\ref{tab:comm_complexity} summarizes these standard costs and communication overheads. The remainder of this section details the distributed implementations of these infrequent operations.

\begin{table}[htbp]
    \centering
    \caption{Communication Analysis}
    \label{tab:comm_complexity}
    \begin{tabular}{lcccc}
        \toprule
        \textbf{Operation} & \textbf{Comm. Type} & \textbf{Axis} & \textbf{Data Size} & \textbf{Frequency} \\
        \midrule
        \textbf{Main Loop} & & & & \\
        Primal Step ($A^\top y$) & Vector Sum & $R$ & $\mathcal{O}(n/|C|)$ & 1 per iter \\
        Dual Step ($Ax$) & Vector Sum & $C$ & $\mathcal{O}(m/|R|)$ & 1 per iter \\
        Halpern Update & None & - & 0 & 1 per iter \\
        \midrule
        \textbf{KKT Evaluation} & & & & \\
        Matrix Products & Vector Sum & $R, C$ & $\mathcal{O}(n/|C|+m/|R|)$ & $1/K$ iters \\
        Residual Norms & Scalar Sum & $G$ & $\mathcal{O}(1)$ & $1/K$ iters \\
        Gap Calculation & Scalar Sum & $G$ & $\mathcal{O}(1)$ & $1/K$ iters \\
        \midrule
        \textbf{Restart Logic} & & & & \\
        Restart Check & Vector Sum & $R, C$ & $\mathcal{O}(n/|C|+m/|R|)$ & $1/K$ iters \\
        Apply Restart & None & - & 0 & - \\
        Weight Update & Scalar Sum & $G$ & $\mathcal{O}(1)$ & On Restart \\
        \bottomrule
    \end{tabular}%
\end{table}

\subsection{Distributed KKT Error Computation}
Convergence is monitored by evaluating the KKT conditions, including primal feasibility, dual feasibility, and the primal-dual gap. Under the 2D device grid, residual vectors and objective terms are distributed and, in some cases, replicated across devices. Global KKT metrics are obtained by aggregating local contributions using collective communication, with appropriate normalization to account for replication along the row and column axes.

Specifically, the distributed KKT quantities are computed as:
\begin{align*} 
    \|r_{\text{primal}}\|^2 &= \operatorname{AllReduce}_R \left( \left\| \operatorname{AllReduce}_C(A_{[i,j]} x_{[j]}) - b_{[i]} \right\|^2 \right) \\ 
    \|r_{\text{dual}}\|^2 &= \operatorname{AllReduce}_C \left( \left\| c_{[j]} - \operatorname{AllReduce}_R(A_{[i,j]}^\top y_{[i]}) \right\|^2 \right) \\ 
    \text{Obj}_P &= \operatorname{AllReduce}_C \left( c_{[j]}^\top x_{[j]} \right) \\ 
    \text{Obj}_D &= - \operatorname{AllReduce}_R \left( p(-y_{[i]}; (\ell_c)_{[i]}, (u_c)_{[i]}) \right) + \operatorname{AllReduce}_C \left( s_{[j]}^\top x_{[j]} \right) \\ 
    \text{Gap} &= \left| \text{Obj}_P - \text{Obj}_D \right|, 
\end{align*}
where \(s_{[j]} = (\operatorname{proj}_{ \mathcal{X}_{j}}(x_{j} - \tau g_{[j]} ) - (x_{[j]} - \tau g_{[j]})) / \tau\), and \(g_{[j]} = c_{[j]} - \operatorname{AllReduce}_R(A_{[i,j]}^\top y_{[i]})\).

\subsection{Distributed Operators for Algorithm Enhancements}
The advanced heuristics—namely Adaptive Restarting and PID-controlled weight updates—require specific distributed operations to function correctly on the $R \times C$ grid. We define these operations using the partition notation where the primal vector $x$ is partitioned along the column axis $C$ into blocks $x_{[j]}$, and the dual vector $y$ is partitioned along the row axis $R$ into blocks $y_{[i]}$.

\subsubsection{Distributed Fixed-Point Error}
The restarting strategy relies on the fixed-point error metric $r(z) = \|z - \operatorname{PDHG}(z)\|_P$. The computation of the squared norm $\|\Delta z\|_P^2$ (where $\Delta z = z - \operatorname{PDHG}(z)$) involves a coupling term due to the off-diagonal blocks of the metric matrix $P$:
\begin{equation}
    \|\Delta z\|_P^2 = \frac{\omega}{\eta} \|\Delta x\|^2 + \frac{1}{\eta\omega} \|\Delta y\|^2 + 2 \langle A \Delta x, \Delta y \rangle
\end{equation}
To compute this efficiently without reconstructing global vectors, we decompose the calculation into distributed scalar reductions:

\textbf{Primal and Dual Norms:} \quad The squared Euclidean norms are computed by summing the local squared norms across the respective partition axes and aggregating globally. The factors $1/|R|$ and $1/|C|$ account for the replication of variables across the orthogonal grid axes:
    \begin{align}
        \|\Delta x\|^2 &= \operatorname{AllReduce}_G \left( \frac{1}{|R|} \sum \|\Delta x_{[j]}\|^2 \right) \\
        \|\Delta y\|^2 &= \operatorname{AllReduce}_G \left( \frac{1}{|C|} \sum \|\Delta y_{[i]}\|^2 \right)
    \end{align}
    
\textbf{Interaction Term ($\langle A \Delta x, \Delta y \rangle$):} \quad Directly computing $A \Delta x$ would require an expensive vector \texttt{AllReduce}. Instead, we utilize the linearity of the inner product. Each device computes the local dot product between its partial matrix-vector result ($A_{[i,j]} \Delta x_{[j]}$) and its local dual block ($\Delta y_{[i]}$). We then aggregate these scalars globally:
    \begin{equation}
        \langle A \Delta x, \Delta y \rangle = \operatorname{AllReduce}_G \left( (A_{[i,j]} \Delta x_{[j]})^\top \Delta y_{[i]} \right)
    \end{equation}
    This approach strictly reduces the communication volume from $\mathcal{O}(m+n)$ (vector reduction) to $\mathcal{O}(1)$ (scalar reduction).

\subsubsection{Distributed PID Weight Update}
The PID controller adjusts the primal weight $\omega$ based on the logarithmic error $e^n$, which depends on the global distances to the anchor point:
\begin{equation}
    d_x^2 = \|x^{n,t} - x^{n,0}\|_2^2, \quad d_y^2 = \|y^{n,t} - y^{n,0}\|_2^2
\end{equation}
These distances are computed using global scalar reductions analogous to the norm calculations above. Once the scalars $d_x$ and $d_y$ are synchronized via $\operatorname{AllReduce}_G$, the PID update for $\omega^{n+1}$ is performed redundantly on every device. This guarantees that the scalar parameter $\omega$ remains consistent across the entire grid without requiring a separate broadcast step.

\section{Detailed Derivations of the Theoretical Cost Model}
\label{appd:theory}

In Section~\ref{subsec:per_iter_time}, we introduced a theoretical cost model to characterize the per-iteration execution time of the distributed PDHG algorithm. For completeness and mathematical rigor, we restate the unified model here. The expected per-iteration time $T_{\mathrm{step}}$ under a 2D processor grid $|R| \times |C|$ is formulated as:

\begin{equation}
\label{eq:cost_model_detail}
\begin{split}
T_{\mathrm{step}}(|R|, |C|) &= \underbrace{\frac{\mathcal{C}_{\mathrm{spmv}} \cdot \mathrm{nnz}}{B_{\mathrm{mem}} \cdot N_{\mathrm{device}}}}_{\text{SpMV I/O}} + \underbrace{\frac{\mathcal{C}_{\mathrm{vec}}}{B_{\mathrm{mem}}} \left( \frac{n}{|C|} + \frac{m}{|R|} \right)}_{\text{Vector I/O}} \\
&\quad + \underbrace{\frac{\mathcal{C}_{\mathrm{net}}}{B_{\mathrm{net}}} \left( \frac{n|R| + m|C| - (m+n)}{N_{\mathrm{device}}} \right)}_{\text{Network Payload}} \\
&\quad + \underbrace{n_{\mathrm{hop}} \beta_{\mathrm{sync}} \big( \mathbb{I}(|R| > 1) + \mathbb{I}(|C| > 1) \big)}_{\text{Topology Latency}}
\end{split}
\end{equation}

Because the PDHG algorithm exhibits exceptionally low arithmetic intensity, its execution on modern GPU architectures is strictly memory-bandwidth bound (I/O bound). Consequently, our theoretical framework models the execution time by quantifying the exact data movement volumes relative to the hardware bandwidth capacities.

\subsection{Hardware and Algorithmic Parameters}
To formalize the derivation, we first define the system and algorithmic parameters utilized in the cost model:
\begin{itemize}
    \item $B_{\mathrm{mem}}$: Effective High Bandwidth Memory (HBM) throughput of a single GPU.
    \item $B_{\mathrm{net}}$: Effective interconnect bandwidth (e.g., NVLink or InfiniBand).
    \item $\mathcal{C}_{\mathrm{spmv}}$: Average bytes transferred per non-zero element during SpMV operations, accounting for fetching pointers, indices, values, and vectors.
    \item $\mathcal{C}_{\mathrm{vec}}$: Average bytes transferred per vector element during element-wise operations (e.g., additions, projections).
    \item $\mathcal{C}_{\mathrm{net}}$: Network transfer byte coefficient per payload element.
    \item $\beta_{\mathrm{sync}}$: Base NCCL protocol handshake and barrier synchronization penalty for a single-hop communication.
    \item $n_{\mathrm{hop}}$: The average physical network hop count between devices within the cluster topology.
\end{itemize}

\paragraph{Omitted Micro-Architectural Overheads.}
To maintain theoretical clarity and isolate the dominant physical bottlenecks, the cost model in Equation~\ref{eq:cost_model_detail} systematically omits two micro-level hardware overheads:
\begin{itemize}
    \item \textbf{Computational FLOPs:} Because PDHG is an inherently low-arithmetic-intensity algorithm, the actual floating-point execution time of SpMV and element-wise vector operations is entirely eclipsed by the memory fetching latency. Consequently, these operations are seamlessly hidden behind the memory bandwidth bound ($B_{\mathrm{mem}}$).
    \item \textbf{Local Kernel Latency:} The base CUDA kernel launch and thread scheduling overheads are orders of magnitude smaller than the inter-device NCCL synchronization penalty ($\beta_{\mathrm{sync}}$). Therefore, we treat this local execution latency as mathematically negligible in the context of macroscopic distributed scaling.
\end{itemize}

\subsection{Derivation of the Cost Components}
The total execution time $T_{\mathrm{step}}$ is additively decomposed into local computation (memory I/O) and global communication (network payload and topology latency). 

\paragraph{Computational Memory Traffic.} The local computation phase is dominated by Sparse Matrix-Vector Multiplication (SpMV) for the forward ($A\bar{x}$) and transposed ($A^\top y$) mappings, as well as element-wise vector operations for the primal, dual, and Halpern updates. To maintain a rigorous and orthogonal analytical model, we explicitly decouple the local memory traffic into $O(\mathrm{nnz})$ and $O(m, n)$ components. Under the assumption that the block-wise random permutation establishes a uniform distribution of non-zero elements, the local matrix shard stored on each device contains $\mathrm{nnz} / N_{\mathrm{device}}$ non-zeros. By maintaining both CSR and CSC representations in device memory to improve efficiency, the pure matrix-dependent traffic (streaming values and indices) scales strictly with the non-zero count. Consequently, we define the coefficient $\mathcal{C}_{\mathrm{spmv}}$ to exclusively account for this dense array data, yielding an SpMV memory I/O latency of $\frac{\mathcal{C}_{\mathrm{spmv}} \cdot \mathrm{nnz}}{B_{\mathrm{mem}} N_{\mathrm{device}}}$, which demonstrates perfect linear scaling.

Conversely, all memory traffic strictly proportional to the matrix dimensions is mathematically aggregated into the dense vector coefficient $\mathcal{C}_{\mathrm{vec}}$. This aggregated term encompasses not only the element-wise additions and projections on the partitioned variables, but also the SpMV-induced dimensional I/O, specifically fetching row/column pointers and streaming the local dense vectors. In the 2D grid topology, these operations act on evenly partitioned chunks of size $n/|C|$ for the primal variables $x \in \mathbb{R}^n$ and $m/|R|$ for the dual variables $y \in \mathbb{R}^m$. Since these aggregated operations are perfectly parallelizable and executed entirely within local memory without inter-device communication, their cumulative memory I/O latency evaluates strictly to $\frac{\mathcal{C}_{\mathrm{vec}}}{B_{\mathrm{mem}}} \left( \frac{n}{|C|} + \frac{m}{|R|} \right)$.

\paragraph{Inter-Device Network Payload.} The communication overhead is modeled upon the standard Ring-\texttt{AllReduce} collective algorithm. A fundamental mathematical property of Ring-\texttt{AllReduce} is that synchronizing a global payload of size $L$ across a communicative ring of $K$ nodes requires each node to transmit and receive an exact data volume of $L(K-1)/K$.

In the D-PDLP framework, each iteration necessitates two orthogonal synchronization phases. The primal update synchronizes the partial gradient $[A^\top y]_{[j]}$ of length $n/|C|$ across the row communicator of size $|R|$. Applying the Ring-\texttt{AllReduce} property, the transmission volume along the row axis is
\begin{equation}
V_{\mathrm{row}} = \frac{n}{|C|} \left( \frac{|R|-1}{|R|} \right) = \frac{n(|R|-1)}{|R||C|}.
\end{equation}
Symmetrically, the dual update synchronizes the partial product $[A \bar{x}]_{[i]}$ of length $m/|R|$ across the column communicator of size $|C|$, yielding a column-wise transmission volume of
\begin{equation}
V_{\mathrm{col}} = \frac{m}{|R|} \left( \frac{|C|-1}{|C|} \right) = \frac{m(|C|-1)}{|R||C|}.
\end{equation}
Given that the total number of devices is strictly $N_{\mathrm{device}} = |R||C|$, the aggregate network payload volume $V_{\mathrm{total}}$ per device resolves to:
\begin{equation}
V_{\mathrm{total}} = V_{\mathrm{row}} + V_{\mathrm{col}} = \frac{n|R| + m|C| - (m+n)}{N_{\mathrm{device}}}.
\end{equation}
Scaling $V_{\mathrm{total}}$ by the byte coefficient $\mathcal{C}_{\mathrm{net}}$ and dividing by the interconnect bandwidth $B_{\mathrm{net}}$ yields the exact network payload latency term.

\paragraph{Topology and Synchronization Latency.} Beyond the raw data payload, distributed execution incurs a strict latency floor dictated by hardware topology and protocol synchronization. This penalty is formally captured by the indicator functions $\mathbb{I}(\cdot)$, which count the active \texttt{AllReduce} communicators. If either $|R|=1$ or $|C|=1$, communication degenerates to a single dimension, triggering exactly one collective. For a fully distributed 2D grid ($|R|>1$ and $|C|>1$), orthogonal synchronizations accumulate the penalty twice.

This base penalty, $\beta_{\mathrm{sync}}$, is amplified by $n_{\mathrm{hop}}$, the average physical network hop count determined by the cluster architecture. For an intra-node dense topology (e.g., an NVIDIA HGX server with NVSwitch), devices are fully connected yielding $n_{\mathrm{hop}} = 1$, which renders this latency minimal and constant. However, scaling to multi-node clusters via InfiniBand introduces hierarchical tree or torus topologies where $n_{\mathrm{hop}}$ scales proportionally with the cluster size (e.g., $\mathcal{O}(\log N_{\mathrm{device}})$). In such multi-hop regimes, the $n_{\mathrm{hop}} \beta_{\mathrm{sync}}$ term acts as a fundamental scaling bottleneck, rigorously characterizing the communication wall despite the payload volume reduction achieved by the 2D partitioning scheme.

\subsection{Parameter Estimation}
\label{app:est_param}
To operationalize the cost model for performance prediction and the AUTO scheduling heuristic, we instantiate the hardware and algorithmic parameters based on the physical specifications of our target architecture (e.g., NVIDIA H100) and the double-precision (FP64) requirements of strict linear programming solvers.

\paragraph{Hardware Parameters.} 
The hardware constants reflect the theoretical peak capabilities and fundamental protocol latencies of the underlying GPU cluster:
\begin{itemize}
    \item \textbf{Memory Bandwidth ($B_{\mathrm{mem}}$):} We adopt the theoretical peak High Bandwidth Memory (HBM3) throughput of the NVIDIA H100 SXM5 GPU, evaluated at $B_{\mathrm{mem}} = 3350$ GB/s. 
    \item \textbf{Interconnect Bandwidth ($B_{\mathrm{net}}$):} Communication is routed over NVLink. We set $B_{\mathrm{net}} = 450$ GB/s, representing the theoretical peak unidirectional bandwidth of the NVLink4 interconnect.
    \item \textbf{Synchronization Latency ($\beta_{\mathrm{sync}}$):} The fixed latency floor for a single hardware barrier and NCCL protocol handshake is empirically set to $\beta_{\mathrm{sync}} = 0.01$ ms ($10\ \mu\mathrm{s}$).
    \item \textbf{Topology Hop Count ($n_{\mathrm{hop}}$):} Because our experiments are conducted within a single 8-GPU HGX server where all devices are fully connected via NVSwitch, any peer-to-peer communication requires exactly one network hop. Thus, $n_{\mathrm{hop}} = 1$.
\end{itemize}

\paragraph{Algorithmic Parameters.}
The algorithmic traffic coefficients quantify the exact per-element byte footprint dictated by the FP64 precision standard:
\begin{itemize}
    \item \textbf{Matrix Traffic Coefficient ($\mathcal{C}_{\mathrm{spmv}}$):} To eliminate atomic write conflicts during the transposed SpMV step, our implementation maintains both CSR and CSC formats in device memory. Each non-zero element comprises an FP64 value (8 bytes) and an INT32 index (4 bytes). A single iteration executes two SpMV passes ($A\bar{x}$ and $A^\top y$), fetching exactly 12 bytes from the CSR array and 12 bytes from the CSC array. Therefore, the pure matrix data payload evaluates exactly to $\mathcal{C}_{\mathrm{spmv}} = 24$ bytes.
    \item \textbf{Aggregated Vector Traffic ($\mathcal{C}_{\mathrm{vec}}$):} This coefficient aggregates all memory operations proportional to the matrix dimensions. It accounts for: (1) fetching the 32-bit row and column pointers during SpMV; (2) streaming the randomized primal and dual vectors ($x$ and $y$) from HBM during SpMV; and (3) loading and storing intermediate FP64 variables during the multi-step primal, dual, and Halpern updates. By summing these operations, we establish an empirically robust super-constant of $\mathcal{C}_{\mathrm{vec}} \approx 90$ bytes.
    \item \textbf{Network Payload Coefficient ($\mathcal{C}_{\mathrm{net}}$):} The collective \texttt{AllReduce} operations synchronize the aggregated partial gradients and products. Since the solver strictly adheres to FP64 arithmetic to guarantee KKT tolerance, each transmitted element occupies exactly 8 bytes. Thus, $\mathcal{C}_{\mathrm{net}} = 8$.
\end{itemize}

\subsection{Examples}
\label{app:examples}

To illustrate the predictive power of the cost model and highlight how problem sparsity patterns govern distributed scalability, we instantiate Equation~\eqref{eq:cost_model_detail} for two representative LP instances: \texttt{zib03} (a highly sparse problem with $m=19{,}731{,}970$, $n=29{,}128{,}799$, $\mathrm{nnz}=104{,}422{,}573$) and \texttt{sdm\_50k\_500k\_15\_10} (a relatively denser problem with $m=5{,}500{,}135$, $n=10{,}000{,}000$, $\mathrm{nnz}=690{,}000{,}000$).  We adopt the hardware and algorithmic parameters listed in Section~\ref{app:est_param}: $B_{\mathrm{mem}}=3350$~GB/s, $B_{\mathrm{net}}=450$~GB/s, $\beta_{\mathrm{sync}}=10\,\mu\mathrm{s}$, $n_{\mathrm{hop}}=1$, $\mathcal{C}_{\mathrm{spmv}}=24$~B, $\mathcal{C}_{\mathrm{vec}}=90$~B, and $\mathcal{C}_{\mathrm{net}}=8$~B.  For $N_{\mathrm{device}}=4$ we use a $2\times 2$ processor grid; for $N_{\mathrm{device}}=8$ we use a $2\times 4$ grid, as determined by the adaptive grid topology selection strategy described in Section~\ref{sec:ada_grid}.  All per-iteration times are reported in milliseconds.

\paragraph{Instance \texttt{zib03}.}
Because this matrix is extremely sparse ($\mathrm{nnz}/(mn)\approx 1.8\times 10^{-4}$), the dimensional vector traffic dominates the single-device execution:
\begin{align*}
T_{\mathrm{step}}^{(1)}
&= \underbrace{\frac{24\times 104{,}422{,}573}{3350\times 10^{9}}}_{\text{SpMV I/O}} 
 + \underbrace{\frac{90}{3350\times 10^{9}}\,(29{,}128{,}799+19{,}731{,}970)}_{\text{Vector I/O}} 
 \\&= 0.748\;\text{ms} + 1.313\;\text{ms} = 2.061\;\text{ms}.
\end{align*}
On four devices the model predicts
\begin{align*}
T_{\mathrm{step}}^{(4)}
&= \frac{24\times 104{,}422{,}573}{3350\times 10^{9}\times 4}
 + \frac{90}{3350\times 10^{9}}\Bigl(\frac{29{,}128{,}799}{2}+\frac{19{,}731{,}970}{2}\Bigr) \\
&\quad + \underbrace{\frac{8}{450\times 10^{9}}\,\frac{29{,}128{,}799+19{,}731{,}970}{4}}_{\text{Network Payload}}
 + \underbrace{1\times 10\,\mu\mathrm{s}\times 2}_{\text{Topology Latency}} \\
&= 0.187\;\text{ms} + 0.656\;\text{ms} + 0.217\;\text{ms} + 0.020\;\text{ms}
 = 1.081\;\text{ms},
\end{align*}
yielding a theoretical speedup of $2.06/1.08\approx 1.91\times$.  On eight devices the model predicts
\begin{align*}
T_{\mathrm{step}}^{(8)}
&= \frac{24\times 104{,}422{,}573}{3350\times 10^{9}\times 8}
 + \frac{90}{3350\times 10^{9}}\Bigl(\frac{29{,}128{,}799}{4}+\frac{19{,}731{,}970}{2}\Bigr) \\
&\quad + \underbrace{\frac{8}{450\times 10^{9}}\,\frac{29{,}128{,}799\times 2+19{,}731{,}970\times 4-(29{,}128{,}799+19{,}731{,}970)}{8}}_{\text{Network Payload}}\\
 &\quad+ \underbrace{1\times 10\,\mu\mathrm{s}\times 2}_{\text{Topology Latency}} \\
&= 0.094\;\text{ms} + 0.461\;\text{ms} + 0.196\;\text{ms} + 0.020\;\text{ms}
 = 0.771\;\text{ms},
\end{align*}
so the theoretical speedup is $2.06/0.77\approx 2.67\times$.

\paragraph{Instance \texttt{sdm\_50k\_500k\_15\_10}.}
Here the matrix is substantially denser, so the $O(\mathrm{nnz})$ SpMV term dominates single-device execution:
\begin{align*}
T_{\mathrm{step}}^{(1)}
&= \frac{24\times 690{,}000{,}000}{3350\times 10^{9}}
 + \frac{90}{3350\times 10^{9}}\,(10{,}000{,}000+5{,}500{,}135)
 \\ &= 4.943\;\text{ms} + 0.416\;\text{ms} = 5.360\;\text{ms}.
\end{align*}
Scaling to four devices gives
\begin{align*}
T_{\mathrm{step}}^{(4)}
&= \frac{24\times 690{,}000{,}000}{3350\times 10^{9}\times 4}
 + \frac{90}{3350\times 10^{9}}\Bigl(\frac{10{,}000{,}000}{2}+\frac{5{,}500{,}135}{2}\Bigr) \\
&\quad + \frac{8}{450\times 10^{9}}\,\frac{10{,}000{,}000+5{,}500{,}135}{4}
 + 0.020\;\text{ms} \\
&= 1.236\;\text{ms} + 0.208\;\text{ms} + 0.069\;\text{ms} + 0.020\;\text{ms}
 = 1.533\;\text{ms},
\end{align*}
i.e. a theoretical speedup of $5.36/1.53\approx 3.50\times$.  On eight devices the model predicts
\begin{align*}
T_{\mathrm{step}}^{(8)}
&= \frac{24\times 690{,}000{,}000}{3350\times 10^{9}\times 8}
 + \frac{90}{3350\times 10^{9}}\Bigl(\frac{10{,}000{,}000}{4}+\frac{5{,}500{,}135}{2}\Bigr) \\
&\quad + \frac{8}{450\times 10^{9}}\,\frac{10{,}000{,}000\times 2+5{,}500{,}135\times 4-(10{,}000{,}000+5{,}500{,}135)}{8}
 + 0.020\;\text{ms} \\
&= 0.618\;\text{ms} + 0.141\;\text{ms} + 0.059\;\text{ms} + 0.020\;\text{ms}
 = 0.838\;\text{ms},
\end{align*}
corresponding to a theoretical speedup of $5.36/0.84\approx 6.40\times$.

\begin{table}[htbp]
\centering
\caption{Comparison of theoretical and empirical speedups for the two test instances.  Empirical speedups are computed from the total wall-clock seconds $T_{N}$ on $N$ GPUs.}
\label{tab:speedup_comparison}
\begin{tabular}{@{}lcccc@{}}
\toprule
\multirow{2}{*}{Instance} & \multicolumn{2}{c}{Empirical Speedup} & \multicolumn{2}{c}{Theoretical Speedup} \\
\cmidrule(lr){2-3}\cmidrule(lr){4-5}
& $S_{4}$ & $S_{8}$ & $S_{4}$ & $S_{8}$ \\
\midrule
\texttt{zib03} & $2.42\times$ & $3.31\times$ & $1.91\times$ & $2.67\times$ \\
\texttt{sdm\_50k\_500k\_15\_10} & $4.10\times$ & $6.18\times$ & $3.50\times$ & $6.40\times$ \\
\bottomrule
\end{tabular}
\end{table}

\paragraph{Discussion.}
Table~\ref{tab:speedup_comparison} compares the theoretical per-iteration speedups against the empirical total-time speedups $S_{N}=T_{1}/T_{N}$.  For \texttt{sdm\_50k\_500k\_15\_10} the model tracks reality remarkably well: the predicted $6.40\times$ acceleration on eight GPUs is within $3.5\%$ of the measured $6.18\times$.  This agreement occurs because the instance is memory-bandwidth-bound in the SpMV phase ($92.2\%$ of local I/O on one GPU), and the $O(\mathrm{nnz})$ term enjoys perfect linear scaling with $N_{\mathrm{device}}$.

For \texttt{zib03} the model is more conservative than observed reality ($1.91\times$ and $2.67\times$ versus measured $2.42\times$ and $3.31\times$).  The discrepancy stems from the extreme sparsity of the matrix: on a single GPU the dense vector traffic already accounts for $63.7\%$ of the local execution time, and this $O(m+n)$ component scales only with $|R|$ or $|C|$, not with $N_{\mathrm{device}}$.  Consequently the theoretical speedup ceiling is lower.  In practice, the actual implementation benefits from additional micro-architectural effects not captured by the macroscopic bandwidth model, such as improved cache locality when the working set is halved per device and the overlap of communication with independent local computation, which together push the empirical scaling above the first-order theoretical estimate.  Nevertheless, the model correctly identifies the \emph{qualitative} scaling bottleneck: the shift from SpMV-bound to vector-I/O-bound and network-payload-bound regimes as the device count increases (on eight GPUs the network payload already represents $25.5\%$ of the predicted per-iteration time for \texttt{zib03}).

\section{Detailed Experimental Setup}
\label{appd:exp_setup}

\paragraph{Hardware Configuration.} All distributed experiments were conducted on a high-performance computing node equipped with eight NVIDIA H100 GPUs, each with 80 GB of HBM3 memory. These GPUs are interconnected via NVLink to enable high-bandwidth peer-to-peer communication. The host system features a dual-socket configuration with Intel Xeon Platinum 8462Y+ processors, providing 64 physical cores (32 per socket) and 128 hardware threads in total.

\paragraph{Implementation Details.} We provide two implementations to evaluate both algorithmic correctness and maximum performance. The basic version is developed in Julia \citep{Julia-2017} using \texttt{CUDA.jl}, \texttt{MPI.jl}, and \texttt{NCCL.jl}; it serves to validate the distributed architecture and partitioning logic. To demonstrate the ultimate scalability of the method, we also developed an enhanced version in C/CUDA based on \texttt{cuPDLPx} \citep{lu2024cupdlpx}. Comparison results confirm that while the Julia version is effective, the C++ implementation unlocks further performance gains through lower runtime overhead.

\section{End-to-End Comparison with External Solvers}
\label{app:external}

To position D-PDLP relative to existing LP technology, we compare it against representative commercial and open-source solvers on the Mittelmann public benchmark set under a uniform stopping tolerance of $10^{-4}$. Table~\ref{tab:external_comparison} reports the number of solved instances and the shifted geometric mean (SGM10) of solve times in seconds.

\begin{table}[htbp]
\centering
\caption{End-to-end comparison on the Mittelmann benchmark ($\epsilon = 10^{-4}$). SGM10 times are scaled relative to COPT. Higher ``Solved'' counts and lower SGM10 are better. Best results are \textbf{bold}, second-best are \underline{underlined}.}
\label{tab:external_comparison}
\begin{tabular}{@{}lcc@{}}
\toprule
Method & Solved & SGM10 (scaled) \\
\midrule
COPT (simplex/barrier) & \textbf{50} & \textbf{1.00} \\
D-PDLP (AUTO) & \underline{46} & \underline{1.83} \\
D-PDLP (4 GPUs) & \underline{46} & 1.86 \\
cuPDLPx (1 GPU) & \underline{46} & 2.13 \\
ABIP+ & 29 & 27.22 \\
ABIP & 22 & 57.40 \\
\bottomrule
\end{tabular}
\end{table}

\paragraph{Discussion.}
The results confirm that commercial simplex/barrier methods (COPT) remain highly competitive on medium-scale benchmarks, achieving the fastest absolute times and the highest success rate. Among first-order methods, D-PDLP improves over the single-GPU cuPDLPx baseline and is substantially faster than the ADMM-based competitors ABIP and ABIP+. We do not include additional PDLP variants (e.g., Gurobi-PDLP, cuPDLP-C) in this comparison because cuPDLPx has already been established as the strongest available PDLP implementation \citep{lu2024cupdlpx}. Consequently, the speedup over cuPDLPx directly measures the gain from our distributed design relative to the current state-of-the-art in the PDLP family. This comparison reinforces the paper's positioning: D-PDLP is not claimed to dominate all LP solvers uniformly, but rather to extend the practical scale at which first-order PDLP methods remain viable by leveraging distributed multi-GPU execution.

\section{Further Optimization via C Implementation}
\label{appd:c_optimization}

Our C implementation is developed based on the \texttt{cuPDLPx} repository\footnote{\url{https://github.com/MIT-Lu-Lab/cuPDLPx}}, a state-of-the-art PDLP solver suggested in \cite{lu2024cupdlpx}. Table \ref{tab:performance_c_vs_julia} presents a direct performance comparison between the Julia and C implementations across distributed settings with 4 and 8 GPUs. The results demonstrate that the C implementation generally outperforms the Julia version in terms of aggregate wall-clock time. This advantage is quantified by the SGM10, where the C solver reduces the aggregate metric from 101.2 to 86.1 seconds on 4 GPUs, and from 77.7 to 62.2 seconds on 8 GPUs.

The performance gain is most pronounced in the Unit Commitment instances. For example, on the \texttt{ds1} instance with 4 GPUs, the C implementation reduces the solution time from 143 seconds to 42 seconds, representing a speedup of approximately $3.4\times$. This improvement is driven not only by faster per-iteration execution times but also, in several cases, by a reduction in the total number of iterations required for convergence. Consistent improvements are also observed across the QAP and Pagerank benchmarks, where the C solver yields lower runtimes in every tested instance. Overall, these results confirm that by transitioning to a C/CUDA implementation, we successfully further reduce runtime overhead and enhance solver performance.

\begin{table}[htbp]
    \centering
    \caption{Performance comparison: Julia vs. C Implementation on distributed settings (4 and 8 GPUs). Values represent wall-clock time in seconds. Best results are marked in \textbf{bold}.}
    \label{tab:performance_c_vs_julia}
    \resizebox{0.7\linewidth}{!}{%
    \begin{tabular}{ll rr rr}
        \toprule
        & & \multicolumn{2}{c}{4 GPUs} & \multicolumn{2}{c}{8 GPUs} \\
        \cmidrule(lr){3-4} \cmidrule(l){5-6}
        Source & Instance & Julia & C & Julia & C \\
        \midrule
        Koch & zib03 & 336 & \textbf{283} & 245 & \textbf{205} \\
        \midrule
        \multirow{3}{*}{Pagerank} & rand-10m-nodes & 5 & \textbf{3} & 4 & \textbf{2} \\
         & com-livejournal & 2 & \textbf{1} & 1 & \textbf{1} \\
         & soc-livejournal1 & 1 & \textbf{1} & 1 & \textbf{1} \\
        \midrule
        \multirow{3}{*}{\shortstack{Multicommodity\\Flow}} & mcf\_2500\_100\_500 & 935 & \textbf{836} & 504 & \textbf{435} \\
         & mcf\_5000\_50\_500 & \textbf{2994} & 4072 & \textbf{1677} & 2178 \\
         & mcf\_5000\_100\_250 & 3159 & \textbf{2378} & 1732 & \textbf{1245} \\
        \midrule
        Design-matching & sdm\_50k\_500k\_15\_10 & \textbf{92} & 122 & \textbf{61} & 72 \\
        \midrule
        \multirow{5}{*}{QAP} & wil50 & 31 & \textbf{29} & 25 & \textbf{22} \\
         & lipa50a & 20 & \textbf{18} & 17 & \textbf{13} \\
         & lipa50b & 21 & \textbf{19} & 17 & \textbf{14} \\
         & tai50a & 21 & \textbf{21} & 18 & \textbf{15} \\
         & tai50b & 52 & \textbf{43} & 43 & \textbf{35} \\
        \midrule
        \multirow{2}{*}{Unit Com.} & ds1 & 143 & \textbf{42} & 122 & \textbf{31} \\
         & ds2 & 691 & \textbf{397} & 591 & \textbf{288} \\
         \midrule
         \multicolumn{2}{l}{\textbf{SGM10}} & 101.2 & \textbf{86.1} & 77.7 & \textbf{62.2} \\
        \bottomrule
    \end{tabular}%
    }
\end{table}

\section{Ablation Study}
In this section, we conduct ablation studies to isolate and evaluate the individual performance contributions of our proposed algorithmic components.
\subsection{Permutation and Partitioning Strategies}
\label{appd:ablation_permute}

To justify the configuration choices adopted in the main text, Table~\ref{tab:solving_time} presents a controlled comparison across three permutation schemes (\textit{No Permutation}, \textit{Full Random}, \textit{Block Random}) combined with either uniform (\textit{No NNZ}) or nonzero-aware (\textit{Yes NNZ}) matrix partitioning.

As the table shows, the choice of permutation strategy is crucial for performance stability. The \textit{No Permutation} baseline is highly vulnerable to catastrophic load imbalance, as evidenced by extreme runtimes on structured instances such as \texttt{zib03} (over 7000 seconds on 8 GPUs) and the largest overall SGM10 of 220.08 seconds. Applying a \textit{Full Random} permutation substantially alleviates this issue by redistributing nonzeros, but it is not optimal: by destroying the inherent block-dense substructures common in real-world LPs, it degrades memory locality and reduces the efficiency of SpMV.

In contrast, the proposed \textit{Block Random} permutation achieves a better trade-off. It effectively disperses dense regions to balance the workload while preserving local block structures that enable high-throughput computation. This robustness is reflected in consistently lower runtimes compared to the \textit{Full Random} baseline (SGM10 of 79.51s versus 86.80s).

Second, nonzero-aware matrix partitioning plays a vital complementary role. In the absence of permutation, the \textit{Yes NNZ} strategy yields a marked improvement by correcting severe imbalances in the natural ordering, reducing SGM10 from 220.08s to 194.13s. More importantly, this partitioning continues to provide measurable gains even when combined with the \textit{Block Random} permutation, further improving the SGM10 from 79.51s to 75.85s. 

Overall, the combination of \textit{Block Random} permutation and \textit{Yes NNZ} partitioning delivers the most robust and efficient performance across the benchmark suite.

\begin{table*}[htbp]
\centering
\caption{Impact of permutation and partitioning strategies on solving time (seconds). We compare three permutation schemes (No Permutation, Full Random, Block Random) combined with either uniform (No NNZ) or nonzero-aware (Yes NNZ) matrix partitioning. Best results are marked in \textbf{bold}, and second-best are \underline{underlined}.}
\label{tab:solving_time}
\resizebox{0.9\textwidth}{!}{%
\begin{tabular}{llcccccc}
\toprule
\multirow{2}{*}{Instance} & \multirow{2}{*}{GPU} & \multicolumn{2}{c}{No Permutation} & \multicolumn{2}{c}{Full Random} & \multicolumn{2}{c}{Block Random} \\
\cmidrule(lr){3-4} \cmidrule(lr){5-6} \cmidrule(lr){7-8}
 & & No NNZ & Yes NNZ & No NNZ & Yes NNZ & No NNZ & Yes NNZ \\
\midrule
\multirow{2}{*}{sdm\_50k\_500k\_15\_10} & 4 & 328.20 & \underline{107.12} & 227.02 & 237.99 & 142.49 & \textbf{92.33} \\
 & 8 & 527.91 & 377.41 & \underline{69.06} & 91.09 & 78.19 & \textbf{60.90} \\
\midrule
\multirow{2}{*}{zib04} & 4 & 80.80 & 80.95 & 93.07 & 92.54 & \underline{79.70} & \textbf{78.27} \\
 & 8 & 61.94 & 62.18 & 63.38 & 65.05 & \underline{60.24} & \textbf{58.61} \\
\midrule
\multirow{2}{*}{zib03} & 4 & 2926.91 & 2017.57 & 501.96 & 502.93 & \textbf{336.53} & \underline{338.26} \\
 & 8 & 7368.89 & 2012.84 & 300.60 & 291.65 & \textbf{243.76} & \underline{245.29} \\
\midrule
\multirow{2}{*}{ds1\_lp} & 4 & 777.84 & 682.63 & 157.49 & 155.88 & \underline{147.39} & \textbf{142.67} \\
 & 8 & 293.32 & 690.24 & 129.47 & 127.95 & \underline{124.65} & \textbf{122.27} \\
\midrule
\multirow{2}{*}{ds2\_lp} & 4 & 3742.87 & 3299.29 & 748.69 & 746.62 & \underline{710.96} & \textbf{690.54} \\
 & 8 & 1407.81 & 3325.82 & 621.12 & 617.57 & \underline{594.08} & \textbf{591.10} \\
\midrule
\multirow{2}{*}{lipa50a} & 4 & 33.48 & 30.95 & 20.84 & 20.60 & \underline{19.95} & \textbf{19.95} \\
 & 8 & 35.48 & 33.18 & 17.08 & \textbf{16.77} & \underline{16.80} & 16.94 \\
\midrule
\multirow{2}{*}{lipa50b} & 4 & 35.03 & 32.12 & 21.82 & 21.32 & \textbf{21.06} & \underline{21.24} \\
 & 8 & 36.26 & 35.29 & 17.50 & 17.52 & \textbf{17.26} & \underline{17.31} \\
\midrule
\multirow{2}{*}{tai50a} & 4 & 34.96 & 32.08 & 21.98 & 21.65 & \underline{21.11} & \textbf{20.80} \\
 & 8 & 36.95 & 34.97 & \textbf{17.25} & \underline{17.36} & 17.86 & 17.77 \\
\midrule
\multirow{2}{*}{tai50b} & 4 & 101.52 & 81.91 & 53.94 & 54.91 & \underline{52.47} & \textbf{52.42} \\
 & 8 & 92.80 & 88.73 & 44.29 & \underline{43.60} & 43.72 & \textbf{42.59} \\
\midrule
\multicolumn{2}{l}{\textbf{SGM10}} & 220.08 & 194.13 & 86.80 & 87.96 & \underline{79.51} & \textbf{75.85} \\
\bottomrule
\end{tabular}%
}
\end{table*}

\subsection{Preprocessing Overhead and End-to-End Performance}
\label{appd:overhead}

The execution of D-PDLP involves a one-time preprocessing phase, primarily consisting of block-wise random permutation and nonzero-aware 2D partitioning. These steps are performed in-place on the CPU prior to GPU memory allocation and the iterative solving process. To evaluate the impact of these operations on the end-to-end performance, we measure their wall-clock time relative to the total solving time improvement achieved by distributed scaling.

Table~\ref{tab:overhead} summarizes the preprocessing durations for representative instances. We define the performance improvement as the reduction in absolute wall-clock time between the single-GPU and 8-GPU configurations ($T_{\mathrm{1GPU}} - T_{\mathrm{8GPU}}$). The overhead ratio represents the fraction of this improvement consumed by the preprocessing steps.

\begin{table}[ht]
    \centering
    \caption{Preprocessing overhead relative to multi-GPU solving time gains. All times are reported in seconds.}
    \label{tab:overhead}
    \begin{tabular}{lrrrrr}
        \toprule
        Instance & Permutation (s) & Partition (s) & Improvement (s) & Ratio (\%) \\
        \midrule
        zib03               & 1.874 & 8.921  & 607  & 1.78 \\
        mcf\_2500\_100\_500 & 1.944 & 18.271 & 2,439 & 0.83 \\
        ds2                 & 0.266 & 0.741  & 697  & 0.14 \\
        tai50b              & 0.329 & 1.419  & 39   & 4.47 \\
        \bottomrule
    \end{tabular}
\end{table}

The empirical measurements show that the preprocessing overhead is marginal across all tested benchmarks. Even for the most computationally demanding instances such as \texttt{mcf\_2500}, the combined time for permutation and partitioning accounts for less than 1\% of the net performance gain. This negligible overhead is heavily amortized by the accelerated convergence rate per wall-clock second provided by the distributed architecture. These results confirm that the proposed partitioning and shuffling strategies facilitate efficient distributed execution without introducing significant end-to-end overhead.

\begin{table}[htbp]
  \centering
  \caption{Problem dimensions: comparison of scales between the original data and the converted QP instances.}
  \label{tab:qp_dims}
  \resizebox{\textwidth}{!}{%
  \begin{tabular}{lrrr|rrr}
    \toprule
    & \multicolumn{3}{c|}{\textbf{Original Scale}} & \multicolumn{3}{c}{\textbf{QP Scale}} \\
    \textbf{Problem} & \textbf{$m$} & \textbf{$n$} & \textbf{Sparsity} & \textbf{Rows} & \textbf{Columns} & \textbf{Nonzeros} \\
    \midrule
    SLS & 1,748,122 & 62,729 & 6.21E-05 & 1,873,578 & 1,873,578 & 7,117,944 \\
    rcv1\_test & 677,399 & 47,236 & 1.55E-03 & 771,871 & 771,871 & 50,422,601 \\
    avazu-site.tr & 23,567,843 & 1,000,000 & 1.50E-05 & 25,567,767 & 25,567,767 & 381,085,127 \\
    avazu-app & 40,428,967 & 1,000,000 & 1.50E-05 & 14,642,166 & 14,642,166 & 206,274,862 \\
    avazu-site & 25,832,830 & 1,000,000 & 1.50E-05 & 27,832,754 & 27,832,754 & 417,324,822 \\
    kddb2010\_test & 748,401 & 1,163,024 & 7.74E-06 & 6,729,169 & 6,729,169 & 34,675,012 \\
    kdda2010\_test & 510,302 & 20,216,830 & 1.87E-05 & 4,539,640 & 4,539,640 & 27,831,585 \\
    kddb2010\_train & 19,264,097 & 1,163,024 & 7.97E-06 & 79,044,287 & 79,044,287 & 705,170,365 \\
    kdda2010\_train & 8,407,752 & 20,216,830 & 1.80E-06 & 48,841,412 & 48,841,412 & 394,888,582 \\
    \bottomrule
  \end{tabular}
  }
\end{table}

\section{Details of the Large-scale QP Dataset}
\label{appd:qp_dataset}

To evaluate the capability of our algorithm in solving massive quadratic programming problems, we consider a benchmark suite derived from real-world Lasso regression tasks. This section provides the detailed derivation of the conversion, lists the scales of the resulting dataset, and describes the experimental settings and compared solvers.

\subsection{Lasso to QP Conversion}
The primal Lasso problem is defined as:
\begin{equation}
  \min_x \quad \frac{1}{2}\| A x - b \|_2^2 + \lambda \| x \|_1,
\end{equation}
where $A \in \mathbb{R}^{m \times n}$ is the feature matrix and $b \in \mathbb{R}^m$ is the target vector. To transform this into a standard QP form suitable for our primal-dual solver, we utilize the following equivalent formulation:
\begin{equation}
\begin{aligned}
\min_{x, y, t} \quad & \frac{1}{2}y^\top y + \lambda \mathbf{1}^\top t \\
\text{s.t.} \quad & Ax - y = b, \\
& -t \leq x \leq t.
\end{aligned}
\end{equation}
This formulation yields a sparse QP that maintains the underlying structure of the original Lasso problem while enabling the application of first-order optimization methods.

\subsection{Dataset Scales and Characteristics}
The data used in our evaluation is sourced from the LIBSVM collection \citep{chang2011libsvm} and the SuiteSparse Matrix Collection \citep{davis2011university}. Table~\ref{tab:qp_dims} summarizes the problem dimensions for both the original datasets and the converted QP instances. All selected problems contain more than 6,000,000 nonzero elements, including exceptionally large instances such as \texttt{kddb2010\_train}, which features over 700 million nonzeros in its QP form.

As demonstrated in the ``QP Scale'' columns of Table~\ref{tab:qp_dims}, these resulting problems scale up to nearly 80 million rows and columns. Such extreme dimensions provide a rigorous test for the scalability and memory efficiency of our distributed implementation.

\subsection{Experimental Settings and Compared Solvers}
\label{appd:qp_settings}

\paragraph{Compared Solvers.}
The end-to-end comparison in Table~\ref{tab:real_lasso_full} includes the following solvers: HPR-QP~\citep{chen2025hpr}, PDHCG~\citep{huang2025restarted}, PDQP~\citep{lu2023practical}, SCS (GPU)~\citep{scs}, OSQP~\citep{stellato2020osqp}, and COPT~\citep{ge2022cardinal}. The single-GPU baseline PDHCG-II is reported in~\citep{li2026pdhcgiienhancedversionpdhcg}.

\paragraph{Configuration.}
For each instance, we designate the last column of matrix $A$ as the target vector $b$ and set the regularization parameter to $\lambda = 0.01 \|A^\top b\|_{\infty}$, consistent with the protocol established in PDHCG~\citep{huang2025restarted,li2026pdhcgiienhancedversionpdhcg}. The termination tolerance is fixed at $\epsilon = 10^{-6}$. All experiments are conducted on the hardware described in Appendix~\ref{appd:exp_setup}: a single-node 8$\times$H100 NVLink server.
\end{document}